\documentclass[12pt]{article}
\usepackage{amssymb,amsmath,amscd,amsthm}
\usepackage{graphicx,epsfig, color,float}

\usepackage{graphicx}
\usepackage[active]{srcltx}

\newtheorem{theorem}{Theorem}[section]

\newtheorem{lemma}[theorem]{Lemma}

\def\be{\color{black}}

\date{}

\begin{document}

\date{}
\title{Metastability in Parabolic Equations and Diffusion Processes with a Small Parameter}
\author{
 M. Freidlin\footnote{Dept of Mathematics, University of Maryland,
College Park, MD 20742, mif@math.umd.edu}, L.
Koralov\footnote{Dept of Mathematics, University of Maryland,
College Park, MD 20742, koralov@math.umd.edu}
} \maketitle

\begin{abstract} 
We study diffusion processes in $\mathbb{R}^d$ that leave invariant a finite collection of manifolds (surfaces or points) in $\mathbb{R}^d$ and small perturbations of such processes. Assuming certain ergodic properties at and near the invariant surfaces, 
we describe the rate at which the process gets attracted to or repelled from the surface, based on the local behavior of the coefficients. 
For processes that include, additionally, a small non-degenerate perturbation, we describe the metastable behavior. Namely, by allowing the time scale to depend on the size of the perturbation, we observe different asymptotic
distributions of the process at different time scales. 

Stated in PDE terms, the results provide the asymptotics, at different time scales, for the solution of the parabolic Cauchy problem when 
the operator that degenerates on a collection of surfaces is perturbed by a small non-degenerate term. This asymptotic behavior switches at a finite number of time scales that are calculated and does not depend on the perturbation. 
\end{abstract}



\section{Introduction} \label{intro}

Consider a diffusion process $X_t$ that satisfies the stochastic differential equation
\[
d X_t = v_0( X_t) dt + \sum_{i=1}^d v_i( X_t) \circ d W^i_t,~~~~~X_0 = x \in \mathbb{R}^d,
\]
where $W^i_t$ are independent Wiener processes and $v_0,...,v_d$ are  bounded $C^3(\mathbb{R}^d)$ vector fields.  
The Stratonovich form  is convenient here since it allows
one to provide a natural coordinate-independent description of the process. The generator of the process is the operator
\begin{equation} \label{opl}
Lu = L_0 + \frac{1}{2} \sum_{i=1}^d L_i^2,
\end{equation}
where $L_i u = \langle v_i, \nabla u \rangle$ is the operator of differentiation along the vector field $v_i$, $i =0,...,d$.  It should
be kept in mind that the process depends on the initial point $x$, yet this is not always reflected in the notation. Sometimes, we will write $X^x_t$ but sometimes, instead, we will use the subscript $x$ to denote the probabilities and expectations  associated with the process starting at $x$. 

Let
$\mathcal{M}_1, ... , \mathcal{M}_m \subset \mathbb{R}^d$ be $C^4$-smooth non-intersecting compact surfaces with dimensions $d_1,...,d_m$, respectively, where $0 \leq d_k < d$, $1 \leq k \leq m$. Points
are considered to be zero-dimensional surfaces. We will assume that each of the surfaces is invariant for the process and that the diffusion restricted to a single surface is an ergodic process (the latter condition is trivially satisfied for $\mathcal{M}_k$ if $d_k = 0$). Denote 
$D =  \mathbb{R}^d \setminus \left( \mathcal{M}_1 \bigcup ... \bigcup \mathcal{M}_m \right)$. We assume that the diffusion matrix is non-degenerate on~$D$.  

For each $k = 1,...,m$  and each $x \in \mathcal{M}_k$, we define $T(x)$ to be the tangent space to $\mathcal{M}_k$ at $x$ (with $T(x)$ being a zero-dimensional space if $d_k = 0$). We
will assume that: 

(a) ${\rm span} (v_0(x),v_1(x),...,v_d(x)) = {\rm span} (v_1(x),...,v_d(x)) = T(x)$ for $x \in \mathcal{M}_1 \bigcup ... \bigcup
 \mathcal{M}_m$;

(b) ${\rm span} (v_1(x),...,v_d(x)) = \mathbb{R}^d$ for $x \in D$. 
\\
This is just a more convenient (and slightly stronger) way to express the assumptions that we already made: (a)~implies that the surfaces are invariant, and the process is ergodic on each surface; 
(b)~means that the diffusion matrix is non-degenerate on $D$. Under Assumption (a), the set $\mathcal{M}_1 \bigcup ... \bigcup
 \mathcal{M}_m$ is inaccessible, i.e., the process $X_t$ starting at $x \in D$
does not it in finite time. 

We will need another assumption, strengthening Assumption (a), concerning the ergodic properties of the process not only on
the surfaces $\mathcal{M}_k$ but in their vicinities. If $\mathcal{M}_k$ is a point, then the assumption concerns the ergodic properties of the
spherical part of the motion near $\mathcal{M}_k$. If $\mathcal{M}_k$ is $(d-1)$-dimensional, then the assumption coincides with Assumption (a). In general, 
the assumption is discussed in detail in the next section. Here, it is sufficient to say that, 
roughly speaking, in a neighborhood of each surface $\mathcal{M}_k$, we can consider an auxiliary process 
obtained from $X_t$ by first 
retaining only the leading terms in the asymptotic expansion of the coefficients at 
$\mathcal{M}_k$ and then taking the projection on the $(d-1)$-dimensional manifold $\mathbf{S}_k = \mathcal{M}_k \times \mathbb{S}^{d-d_k-1}$
(with $\mathbb{S}^{d-d_k-1}$ denoting the $(d-d_k-1)$-dimensional sphere) obtained by factoring out the distance to the surface $\mathcal{M}_k$. The resulting process needs to be ergodic on $\mathbf{S}_k$. 
Under such an assumption,  
each surface will be classified
as either attracting or repelling for $X_t$, depending (roughly speaking) on whether 
$\mathrm{P}_x (\lim_{t \rightarrow \infty} {\rm dist}(X_t, \mathcal{M}_k) = 0) > 0$  for each $x$ sufficiently close to $\mathcal{M}_k$ (attracting surface) or 
$\mathrm{P}_x (\lim_{t \rightarrow \infty} {\rm dist}(X_t, \mathcal{M}_k) = 0) = 0$  for each $x \notin \mathcal{M}_k$ (repelling surface). 

Let $U_k$, $1 \leq k \leq m$,  be open sets with disjoint closures such that $\mathcal{M}_k \subset  U_k$, $1 \leq k \leq m$. With each attracting surface $\mathcal{M}_k$, we will associate a number $\gamma_k >0$ such
that the probability that the process $X_t$ starting at a distance $r > 0$ from $\mathcal{M}_k$ ever exits $U_k$ is of order 
$r^{\gamma_k}$ when $r = {\rm dist}(x, \mathcal{M}_k) \downarrow 0$. 
If $\mathcal{M}_k$ is repelling, we can associate a number ${\gamma}_k < 0$ to it
such that the probability that the process $X_t$ starting at a fixed point $x \in U_k \setminus \mathcal{M}_k$  reaches the $r$-neighborhood of $\mathcal{M}_k$ before exiting $U_k$ is of order $r^{ -{\gamma}_k}$ when $r \downarrow 0$.  
 The numbers $\gamma_k$ can be
found by solving a certain non-linear spectral problem that involves an operator obtained from $L$ by linearizing 
the coefficients near $\mathcal{M}_k$. In order to avoid unnecessary technical details, we will assume that $\gamma_k \neq 0$. 

 Next, we perturb the process $X_t$ by a small non-degenerate diffusion. 
The resulting process $X^{\varepsilon}_t$ satisfies
the following stochastic differential equation:
\[
d X^{\varepsilon}_t = (v_0 + \varepsilon^2   \tilde{v}_0) (X^{\varepsilon}_t ) dt + \sum_{i=1}^d v_i(X^{\varepsilon}_t ) \circ d W^i_t + \varepsilon 
\sum_{i =1}^d \tilde{v}_i (X^{\varepsilon}_t ) \circ d \tilde{W}^i_t,~~~X^{\varepsilon}_t = x \in \mathbb{R}^d,
\]
where $\tilde{W}^i_t$ are independent Wiener processes (also independent of all $W^i_t$), and $\tilde{v}_0,...,\tilde{v}_d$ are 
bounded 
$C^3(\mathbb{R}^d)$ vector fields.   In order to make our assumption on the 
non-degeneracy of the
perturbation more precise, we state it as follows: 

(c) ${\rm span} (\tilde{v}_1(x),...,\tilde{v}_d(x)) = \mathbb{R}^d$ for $x \in \mathbb{R}^d$. 

We also include a condition that ensures that the processes $X_t$ and $X^{\varepsilon}_t$ (with sufficiently small $\varepsilon$) are positive-recurrent. Namely, we assume that:

(d)  $v_0,...,v_d, \tilde{v}_0,...,\tilde{v}_d$ are bounded and $\langle v_0(x), x \rangle < -c \|x\|$ for some $c > 0$ and all sufficiently large
$\|x\|$.  
\\

 The generator of the process $X^{\varepsilon}_t$ is the operator $L^\varepsilon = L + \varepsilon^2 \tilde{L}$ with
\[
\tilde{L}u = \tilde{L}_0 + \frac{1}{2} \sum_{i=1}^d \tilde{L}_i^2,
\]
where $\tilde{L}_i u = \langle \tilde{v}_i, \nabla u \rangle$ is the operator of differentiation along the vector field $\tilde{v}_i$, $i =0,...,d$. 

The process $X_t$  starting at $x$ is approximated well by $X^\varepsilon_t$  (with the same starting point)  on finite time intervals as $\varepsilon \downarrow 0$. However, due to the presence of the small non-degenerate component,  
$X^{\varepsilon}_t$ cannot get forever trapped in a vicinity of an attracting surface. 
%
Our goal  is to describe the behavior of the  processes $X^{\varepsilon}_t$ 
at 
times $t(\varepsilon)$ that tend to infinity as $\varepsilon \downarrow 0$.
Roughly speaking, one can divide the quadrant $(0,\infty) \times (0, \infty)$ into a finite number
of domains such that $X^{\varepsilon}_{t(\varepsilon)}$ has a limiting distribution (that depends on the initial point)
when $(1/\varepsilon, t(\varepsilon))$ approaches infinity without leaving a given domain.  These are the metastable distributions for the process. 

 The metastable distributions will be seen to be linear combinations of certain probability measures, Namely, let 
$\pi_k$ be the invariant measure of the unperturbed process $X_t$ restricted to $\mathcal{M}_k$. Such a measure exists and is unique due to
Assumption (a). If at least one of the surfaces $\mathcal{M}_k$ is attracting, then the metastable distributions are linear combinations of
those $\pi_k$ that correspond to the attracting surfaces. If all the
surfaces $\mathcal{M}$ are repelling, then the distribution of $X^{\varepsilon}_{t(\varepsilon)}$ converges, as $\varepsilon \downarrow 0$ and $t(\varepsilon) \rightarrow \infty$ (irrespectively of how slowly $t(\varepsilon)$ goes to infinity,  provided that the initial point belongs to $D$),  to the invariant measure $\mu$ for the unperturbed process on $D$ (this measure exists and is unique, as shown in Section \ref{setwo}). The precise formulation of the result on metastability and the formulas for the coefficients corresponding to $\pi_k$ at each time scale will be provided in Section~\ref{mdis}. 

Let us contrast the differences in the formation of metastable states in the case of random perturbations of dynamical systems 
(see \cite{FW}, Ch.6) and the case of random perturbations of degenerate diffusions considered in this paper. In the former case, the switching times have exponential asymptotics and depend on the initial point, more precisely, on the asymptotically stable state of the 
unperturbed system attracting the initial point. The metastable states are described by a hierarchy of cycles defined by the quasi-potentials.
In the case under consideration here, there is no hierarchy of cycles (although the construction of the hierarchy of Markov chains 
from \cite{FK-MC} formally applies if $D$ is not connected), the switching times have
power asymptotics and are the same for all initial points (if $D$ is connected), so only the weights with which different 
measures $\pi_k$ are included in a metastable distribution depend on the initial point. 

Let us briefly discuss the PDE interpretation of the results on metastability.
Consider the parabolic Cauchy problem: 
\begin{equation} 
\begin{split} \label{direq2} 
\frac{\partial u^\varepsilon (t,x)}{\partial t}&  = L^\varepsilon u^\varepsilon (t, x),~~t > 0, x \in \mathbb{R}^d; \\
u^\varepsilon(0, x) & =   
g(x),~~x \in \mathbb{R}^d,
\end{split}
\end{equation}
where $g$ is a bounded continuous function. The solution $u^\varepsilon$ of problem (\ref{direq2})  can be written as an  expectation
of a functional of the corresponding diffusion processes that depends on the parameter $\varepsilon$, and the dependence of the limit of  
$u^\varepsilon(t(\varepsilon),x)$  on the asymptotics of $ t(\varepsilon)$ is a manifestation
of metastability for the underlying diffusion.  Without loss of generality, we can assume that the surfaces
$\mathcal{M}_1,..., \mathcal{M}_m$ are labeled in such an order that $\gamma_1 \geq  ... \geq  \gamma_m$. Let $\overline{m}$ be such 
that $\gamma_{\overline{m}} > 0 > \gamma_{\overline{m}+1}$ (we put $\overline{m}  = 0$ if $\gamma_1 < 0$ and $\overline{m} = m$ if $\gamma_m 
> 0$). 

If a manifold $\mathcal{M}_k \subset \mathbb{R}^d$ has dimension $d_k < d-1$, then we should make an assumption concerning non-degeneracy of operator $L$ in the vicinity of $\mathcal{M}_k$. Such an assumption (assumption (e)), strengthening assumption (a), is introduced in the next section.
 One of the main results of this paper, stated in PDE terms, is the following.
\begin{theorem} \label{mnt2}  Let all invariant manifolds have dimensions $d_k < d-1$. (See the discussion below for the case $d_k = d-1$.) Under assumptions (a)-(e) on the coefficients of the operator $L^\varepsilon$,
the asymptotics of solutions to the initial-boundary value problem (\ref{direq2}) is as follows. 

(a) If at least one surface is attracting (i.e., $\gamma_1 > 0$), then there exist 
 $p^{x,l}_k \geq 0$,   $1 \leq l \leq \overline{m}$, $1 \leq k \leq l$,  that are continuous 
functions of $x \in \mathbb{R}^d$, with $\sum_{k =1}^{l}p^{x,l}_k = 1$, such that, for $x \in D$,  
 \[
	\lim_{\varepsilon \downarrow 0} u^\varepsilon(t(\varepsilon), x) = \sum_{k = 1}^{\overline{m}} p^{x,\overline{m}}_k 
	\int_{\mathcal{M}_k} g  d \pi_k,~~
	if~~~1 \ll t(\varepsilon) \ll \varepsilon^{-\gamma_{\overline{m}}},
	\]
  \[
	\lim_{\varepsilon \downarrow 0} u^\varepsilon(t(\varepsilon), x) = \sum_{k = 1}^{l} p^{x,l}_k \int_{\mathcal{M}_k} g  d \pi_k,~~
	if~~~\varepsilon^{-\gamma_{l+1}} \ll t(\varepsilon) \ll \varepsilon^{-\gamma_{l}}~~with~~1 \leq l < \overline{m},
	\]
  \[
	\lim_{\varepsilon \downarrow 0} u^\varepsilon(t(\varepsilon), x) =  \int_{\mathcal{M}_1} g  d \pi_1,~~
	if~~~ t(\varepsilon) \gg \varepsilon^{-\gamma_{1}}~~ and ~~\gamma_1 > \gamma_2. 
	\]

(b) If all the components of the boundary are repelling (i.e., $\gamma_1 < 0$), then,  for $x \in D$,  
\[
	\lim_{\varepsilon \downarrow 0} u^\varepsilon(t(\varepsilon), x) =  \int_{D} g  d \mu,~~
	if~~~ t(\varepsilon) \gg 1.
\]

If, additionally, we assume that $t(\varepsilon) \gg |\ln(\varepsilon)|$ (which is relevant for the first limit in (a) and the limit in (b)), then all the above limits hold for all $x \in \mathbb{R}^d$. 

\end{theorem}

Observe that the measures $\pi_1,...,\pi_m$ and $\mu$ do not depend on the coefficients of the perturbation $\tilde{L}$. Neither do 
the coefficients $p^{x,l}_k$ in the representation of the solution, as will follow from the arguments below. Thus the asymptotic
behavior of 
solutions to the Cauchy problem (\ref{direq2}) (for the time scales satisfying
the assumptions of Theorem~\ref{mnt2}) is determined exclusively by the unperturbed operator $L$ and does not depend on the 
perturbation $\tilde{L}$, assuming that the latter satisfies the above assumptions.  

The question on the asymptotic behavior of the solutions to the Cauchy problem is closely related to the study of parabolic initial-boundary value problems and elliptic problems with degeneration on the boundary. Certain classes of PDEs with degeneration on the boundary (as a rule, without a perturbation) were considered in 
\cite{Fich}, \cite{Has}, \cite{OR}, \cite{F85}  (see also references therein). The boundary of the domain is usually assumed to be a $(d-1)$-dimensional surface or a union of such surfaces. In a paper of M. Pinsky (\cite{MP}), a criterion for stochastic stability (similar to our notion of an attracting surface) was given for processes that degenerate in a vicinity of a point and certain sufficient conditions were provided for a surface to be attracting.  

Processes with degeneration on hypersurfaces and the asymptotics of solutions to the corresponding PDEs were considered in our recent papers \cite{FK21} and \cite{FKn2}. In \cite{FK21}, we described the limiting behavior, as the size of the perturbation goes to zero, for solutions to the elliptic Dirichlet problem with degeneration on the boundary, and in \cite{FKn2} we proved a result similar to 
Theorem~\ref{mnt2} for the initial-boundary value problems with the Dirichlet and Neumann boundary data. The difference between the latter results and the current paper is that now the  invariant surfaces have dimension lower than  $(d-1)$. Thus the unperturbed process may ``wind around" the surface, requiring a more delicate analysis of the local behavior. 
Another complication is that the perturbed process, even though it is non-degenerate, does not hit a surface whose dimension is smaller than
$(d-1)$, and thus we can't conveniently consider such stopping times in order to analyze various transition probabilities and transition times.  Some of the lemmas proved in the current paper concerning the behavior of the unperturbed and perturbed processes near an invariant surface remain valid  in the case when the surface is $(d-1)$-dimensional. However, if some of the repelling surfaces are $(d-1)$-dimensional, then the domain $D$ is not connected, and the unperturbed process may have multiple invariant measures on $D$. The limiting invariant measure for the perturbed process will be a linear combination of such measures. A description of the metastable behavior in such a case will be the subject of a forthcoming paper.   

Let us briefly discuss the structure of the paper. In Section~\ref{sse1a}, we introduce a change of variables near an invariant surface and prove a crucial lemma  (Lemma~\ref{spec}) on a non-linear spectral problem related to the generator of the process.  The lemma defines the scaling exponents 
$\gamma_k$ that later find their way into the formulas for the transition times and transition probabilities near the invariant surfaces.  
These new variables and the lemma
later allow us to approximate certain functions of the processes $X_t$ and $X^\varepsilon_t$ by martingales in the vicinities of the 
surfaces $\mathcal{M}_k$ and thus study transition probabilities and transition times in the vicinity of $\mathcal{M}_k$. Transition probabilities for the unperturbed and perturbed processes are studied in Sections~\ref{tpni} and \ref{tpni2}, respectively. Transition times
near an invariant surface are studied in Section~\ref{sttim}. In Section~\ref{dpen}, we prove that the distribution of the unperturbed or perturbed process, as it leaves a fixed neighborhood of an attracting surface, does not substantially depend on the starting point if the starting point is chosen asymptotically close to the surface. This observation is important for the description of the long-time behavior. 
In Section~\ref{mdis}, we combine the above results to describe the behavior of the process $X^\varepsilon_t$ at different time scales, thus 
providing a probabilistic version (Theorems~\ref{iuo}, \ref{btto}, \ref{ltsc}, and \ref{hujuh4}) of Theorem~\ref{mnt2}. 
The long time behavior of the unperturbed process $X_t$ is described earlier, in Section~\ref{ldup}.    

\section{Asymptotic behavior of the unperturbed processes} \label{setwo}

\subsection{Structure of the generator near an invariant surface} \label{sse1a}

In this section, we consider the behavior of the  process $X_t$ near a single surface $\mathcal{M}_k$. Since the surface is fixed, we denote it by $\mathcal{M}$, thus dropping the subscript, and denote its dimension by $d' = d_k$.  
We assume that $d - d' > 1$ (the situation when $d - d' = 1$ has been largely considered in \cite{FKn2}, and the current results can also be seen to hold in this case, but some of the notation in the proofs would need to be modified). 

For $x \in \mathcal{M}$, let $\Pi_x$ be the $(d - d')$-dimensional affine space in $\mathbb{R}^d$ that  is orthogonal to $\mathcal{M}$ 
at $x$. Let $\mathbb{S}$ denote the $(d-d'-1)$-dimensional unit sphere and $\mathbf{S} =   \mathcal{M} \times \mathbb{S}$ be the manifold of dimension $d-1$ that is the direct product of  our surface $\mathcal{M}$ and the sphere. 
Note that, in the case when $\mathcal{M}$ is a single point, $\Pi_x$ is the entire space $\mathbb{R}^d$ and $\mathbf{S} = \mathbb{S}$.

For each $x$ in a sufficiently small
neighborhood of $\mathcal{M}$, there are unique $m(x) \in \mathcal{M}$ and $z(x) \geq 0$ such that $x = m(x) + z(x) n(x)$, where $n(x)$ is a vector of unit length with $m(x)$ as a starting point and an end point in $\Pi_{m(x)}$ (that is $n(x)$ is orthogonal to $\mathcal{M}$ at 
$m(x)$).

 Let $\varphi(x) = (m(x), n(x), z(x))$; for all sufficiently small $r > 0$, this is a bijection between an $r$-neighborhood of $\mathcal{M}$, denoted by $\mathcal{M}^r$, and the set $\mathcal{M} \times \mathbb{S} \times [0, r) = \mathbf{S} \times [0,r)$. We will write 
$(m,n,z) = (y,z) = \varphi(x) \in \mathbf{S} \times [0,r)$ to denote the representation of $x$ in this new set of coordinates. 

Next, let us represent the vector fields $v_0,...,v_d$ in $(y,z)$ coordinates.  For $x \in \mathcal{M}^r \setminus \mathcal{M}$, the vector field $v_i(x)$, when 
viewed in $(y,z)$ coordinates, is a tangent vector to the manifold $\mathbf{S} \times (0,r)$.  For $y = (m,n) \in \mathbf{S}$, the tangent space $T^{\mathbf{S}}_y$ can be represented as $T^{\mathbf{S}}_y = T^\mathcal{M}_m \times T^{\mathbb{S}}_n$. Thus $v_i$ can be represented, in the coordinate form, as
\[
v_i (y,z) = v_i(m,n,z) =  (v^1_i(m,n,z), v^2_i(m,n,z), v^3_i(m,n,z)),
\]
where 
\[
v^1_i(m,n,z)\in  T^\mathcal{M}_m,~ v^2_i(m,n,z)\in T^{\mathbb{S}}_n,~v^3_i(m,n,z) \in \mathbb{R}.
\]
Using the Taylor expansion of $v_i(x)$ around $m(x)$, we can find a square $(d-d') \times (d-d')$ matrix $M_i(m)$ that depends smoothly on $m$ such that,
 for $x \in \mathcal{M}^r \setminus \mathcal{M}$,
\[
v^1_i(m,n,z) = v_i(m) + z g_{i,1}(m,n,z),
\]
\[
v^2_i(m,n,z) = M_i(m) n -  \langle M_i(m) n, n \rangle n + z g_{i,2}(m,n,z),
\]
\[
v^3_i(m,n,z) = z \langle M_i(m) n, n \rangle + z^2 g_{i,3}(m,n,z),
\]
where $g_{i,1}$, $g_{i,2}$, $g_{i,3}$ are $C^2$ functions with bounded derivatives on $\mathbf{S} \times (0,r)$. 
For $y \in \mathbf{S}$, let $w_i(y) = ( v_i(m), M_i(m) n -  \langle M_i(m) n, n \rangle n) \in T^\mathcal{M}_m \times T^{\mathbb{S}}_n = 
T^{\mathbf{S}}_y $. 
 Define
the differential operator $L_y$ on $\mathbf{S}$ via 
\begin{equation} \label{opl22}
L_y = \overline{L}_0 + \frac{1}{2} \sum_{i=1}^d \overline{L}_i^2,
\end{equation}
where $\overline{L}_i $ is the operator of differentiation along the vector field ${ {w}}_i$, $i =0,...,d$.
The subscript $y$ is used here to stress that the operator acts only in the $y$ variables.  Intuitively, this operator 
is the generator for the leading part of the motion near the surface, where the distance from the surface is disregarded.  We will make the following assumption:

(e) $L_y$ is elliptic on $\mathbf{S}$, i.e., ${\rm span}({ {w}}_1,...,{ {w}}_d) = T^\mathbf{S}_y$ for
each $y \in \mathbf{S}$.

Note that the operator
$L_y$ acting in the $y$ variables can be applied to functions of $(y,z)$ by treating $z$ as a parameter. The following lemma provides a
representation for the operator $L$ in $(y,z)$ coordinates in a unit neighborhood of the origin, while specifying the main terms as $z \downarrow 0$.  
\begin{lemma} \label{geneNN} The generator of the process $X_t$ in $(y,z)$ coordinates can be written as:
\begin{equation} \label{oprww}
L u  = L_y u  +  z^2 \alpha(y) \frac{\partial^2 u}{\partial z^2} +   z  \beta(y) \frac{\partial u}{\partial z}  + z {\mathcal{D}}_y \frac{\partial u}{\partial z} + R u 
\end{equation}
with
\[
R u = z {\mathcal{K}}_y u   + z^2 \mathcal{N}_y \frac{\partial u}{\partial z}  +  z^3 \sigma(y,z) \frac{\partial^2 u}{\partial z^2},
\]
where ${\mathcal{D}}_y$ is a  differential operator with first-order derivatives in $y$, whose coefficients depend only on the $y$ variables, ${\mathcal{K}}_y$ is a  differential operator  on $\mathbf{S} \times (0,1)$ with first- and second-order 
derivatives in $y$,  ${\mathcal{N}}_y$ is a differential operator on 
$\mathbf{S} \times (0,1)$ with first-order derivatives in $y$ and a potential term. All the 
operators have  coefficients that are $C^1$ and bounded on their respective domains,  
while $\alpha, \beta \in C^1(\mathbf{S})$, $\sigma \in C^1(\mathbf{S} \times (0,1))$ and $\sigma$ is  bounded. 
\end{lemma}
\proof 
Each term of the operator $L$ defined in (\ref{opl}) can
be considered separately, so let us first consider the case when
\[
Lu = \frac{1}{2} \frac{\partial}{\partial v}(\frac{\partial u}{\partial v}),
\] 
where $v$ is one of the vector fields $v_1,...,v_d$. 
Let $v(x) =  a(x) + b(x)$, where, in $(m,n, z)$ coordinates, $a$ and $b$ are represented as
\begin{equation} \label{vfa}
a(m,n,z) =  ( v(m), M(m) n -  \langle M(m) n, n \rangle n, z \langle M(m) n, n \rangle) \in T^\mathcal{M}_m \times T^{\mathbb{S}}_n \times
\mathbb{R},
\end{equation}
\begin{equation} \label{vfb}
b(m,n,z) = (z g_1(m,n,z), z g_2(m,n,z), z^2 g_3(m,n,z)) \in T^\mathcal{M}_m \times T^{\mathbb{S}}_n \times
\mathbb{R},
\end{equation}
where $M$ is a square matrix and $g_1$, $g_2$, $g_3$ are $C^2$ functions with bounded derivatives on $\mathbf{S} \times (0,r)$.
Thus
\begin{equation} \label{lu11}
Lu = 
\frac{1}{2} ( \frac{\partial}{\partial a}(\frac{\partial u}{\partial a})+   \frac{\partial}{\partial a}(\frac{\partial u}{\partial b}) +
 \frac{\partial}{\partial b}(\frac{\partial u}{\partial a}) +  \frac{\partial}{\partial b}(\frac{\partial u}{\partial b})).
\end{equation}
Let $\overline{L}$ be the operator of differentiation along the vector field 
 $w(y) = ( v(m), M(m) n -  \langle M(m) n, n \rangle n) \in T^\mathcal{M}_m \times T^{\mathbb{S}}_n = 
T^{\mathbf{S}}_y $ on $\mathbf{S}$.
By (\ref{vfa}),
\[
\frac{\partial}{\partial a}(\frac{\partial u}{\partial a}) =
\]
\[
\overline{L}^2 u + z^2
\langle M(m) n, n \rangle^2 \frac{\partial^2 u}{\partial z^2} +
z (\overline{L} ( \langle M(m) n, n \rangle)   +  \langle M(m) n, n \rangle^2 ) \frac{\partial u}{\partial z} +
 2 z  \langle M(m) n, n \rangle \overline{L} (\frac{\partial u}{\partial z}).
\]
As follows from (\ref{vfa}), (\ref{vfb}), the last three terms in (\ref{lu11}) satisfy the conditions imposed on the remainder term $R$
in (\ref{oprww}).

Now consider the case when the operator $L$ defined in (\ref{opl}) is of the first order, i.e., 
\[
Lu =   \frac{\partial u}{\partial v},
\] 
where $v = v_0$.  Thus
\[
Lu = \frac{\partial u}{\partial a}  + \frac{\partial u}{\partial b}.
\]
By (\ref{vfa}),
\[
 \frac{\partial u}{\partial a} = \overline{L} u + z  \langle M(m) n, n \rangle \frac{\partial u}{\partial z},
\]
while ${\partial u}/{\partial b}$ satisfies the conditions imposed on the remainder term $R$. 

Thus the operator $L$ given in (\ref{opl}) (and containing $d+1$ terms) can be written in the form (\ref{oprww}) with 
\[
\alpha(y) = \frac{1}{2} \sum_{i=1}^d \langle M_i(m) n, n \rangle^2,~~ 
\]
\[
\beta(y) = \langle M_0(m) n, n \rangle + \frac{1}{2}  \sum_{i=1}^d (\overline{L}_i ( \langle M_i(m) n, n \rangle)   +  \langle M_i(m) n, n \rangle^2 ),
\]
\[
{\mathcal{D}}_y = \sum_{i=1}^d  \langle M_i(m) n, n \rangle \overline{L}_i.
\]
\qed
\\
{\bf Remark.} Lemma~\ref{geneNN} shows that the operator $L$  is approximated, for small $z$, by the operator
\begin{equation} \label{kok}
K u  = L_y u  +  z^2 \alpha(y) \frac{\partial^2 u}{\partial z^2} +   z  \beta(y) \frac{\partial u}{\partial z}  + 
z {\mathcal{D}}_y \frac{\partial u}{\partial z}.
\end{equation}
The latter operator has a useful homogeneity property: $(Ku)(y,cz) = K (u(y,cz))$ for each $y$ and $c >0$. 
\\

Let $\pi$ be the probability measure on $\mathbf{S}$ that is invariant for the process governed by the operator $L_y$. By Assumption (e), such measure exists and is unique.  Let $\mathbf{\pi}_{\mathcal{M}}$ be the marginal distribution of $\pi$ on $\mathcal{M}$. 
Let $f \in C^2(\mathcal{M})$. Such a function $f$ can be extended to $\mathbf{S}$ by defining  
$\hat{f}(m,n) = f(m)$. Observe that, by (\ref{opl22}), $L_y \hat{f} (m,n)= L f(m)$ for each $(m,n) \in \mathbf{S}$, where $L$ is viewed as a restriction of the operator in (\ref{opl}) to $\mathcal{M}$ 
(equivalently, as the generator of the process $X_t$  on the invariant surface $\mathcal{M}$). Therefore,
\[
0 = \int_{\mathbf{S}} L_y \hat{f} d \pi = 
\int_{\mathcal{M}} L f d 
\pi_{\mathcal{M}}.
\]
This shows that $\pi_{\mathcal{M}}$ seves as the invariant measure for the process $X_t$  on the invariant surface $\mathcal{M}$.

Define 
\[
\bar{\alpha} = \int_{\mathbf{S}} \alpha(y) d\pi(y),~~\bar{\beta} = \int_{\mathbf{S}} \beta(y) d\pi(y) .
\]
 We will see that  if $\bar{\alpha} > \bar{\beta} $,
then $\mathrm{P}_x (\lim_{t \rightarrow \infty} {\rm dist}(X_t, \mathcal{M}) = 0) > 0$ 
for each $x$ in a sufficiently small neighborhood of $\mathcal{M}$. 
If $\bar{\alpha} < \bar{\beta}$, then this
probability is zero unless $x \in \mathcal{M}$. We will assume that $\bar{\alpha} \neq \bar{\beta}$   and will
refer to $\mathcal{M}$ as attracting if $\bar{\alpha} > \bar{\beta}$ and repelling if $\bar{\alpha} < \bar{\beta}$.

The next lemma will be instrumental in analyzing the transition probabilities for the process $X_t$ as well as the transition times 
for the process $X^{ \varepsilon}_t$ in the vicinity of $\mathcal{M}$. It was first proved in \cite{FKn2}. Since the proof is short, we 
include it here for the sake of completeness. 

\begin{lemma} \label{spec}
If $\bar{\alpha} > \bar{\beta}$ ($\bar{\alpha} < \bar{\beta}$), then there exist a constant $\gamma >0$ ($\gamma <0$) and a positive-valued function 
$\varphi  \in C^2(\mathbf{S})$ satisfying $\int_\mathbf{S} \varphi d \pi = 1$ such that
\begin{equation} \label{mg}
L_y \varphi (y) + \gamma (\gamma -1) \alpha (y) \varphi (y) + \gamma \beta (y)  \varphi (y) + \gamma \mathcal{D}_y \varphi (y) = 0,~~~y \in
 \mathbf{S}.
\end{equation}
Such $\gamma$ are $\varphi$ are determined uniquely. 
\end{lemma}
\proof Let $\lambda_\gamma$ be the top eigenvalue for $M(\gamma)$, where $M(\gamma)$ is the operator in the left hand side of (\ref{mg}).
Since $L_y + \gamma \mathcal{D}_y$ is an elliptic operator on a compact manifold, this eigenvalue is simple, as follows from the 
Perron-Frobenius Theorem by considering the corresponding parabolic semigroup $T_t$, $t \geq 0$,  with the kernel of $T_t$ positive for each $t > 0$. 
Therefore  $\lambda_\gamma$ depends smoothly on the parameter $\gamma$  (\cite{Kato}). Moreover, the 
 corresponding eigenfunction $\varphi_\gamma$ can be chosen so that it depends smoothly on $\gamma$ and $\varphi_0 \equiv 1$. 
Differentiating the equality
$ M(\gamma) \varphi_\gamma = \lambda_\gamma \varphi_\gamma $ in $\gamma$, we obtain
\[
M(\gamma) \varphi'_\gamma + ((2 \gamma - 1) \alpha  + \beta +  \mathcal{D}_y) \varphi_\gamma = \lambda'_\gamma \varphi_\gamma +  \lambda_\gamma \varphi'_\gamma.
\]
Put $\gamma = 0$ and integrate both sides with respect to the measure $\pi$. 
Note that $\int_\mathbf{S} \varphi_0 d \pi = 1$, $\mathcal{D}_y \varphi_0 = 0$, and $\lambda_0 = 0$.
Since $M(0) = L_y$ and $\pi$ is invariant for the process generated by $L_y$, we have
$\int_\mathbf{S} M(0) \varphi'_0 d \pi  = 0$. Thus 
\[
 \lambda'_0  = \int_\mathbf{S} (-\alpha + \beta) d \pi = \bar{\beta} - \bar{\alpha}.
\]
Let us assume that $\bar{\beta} <  \bar{\alpha}$ (the case when $\bar{\beta} >  \bar{\alpha}$ can be handled similarly). 
Then $ \lambda'_0 < 0$. Observe that $\lambda_0 = 0$. 
Let us demonstrate that
  $\lim_{\gamma \rightarrow \infty} \lambda_\gamma = +\infty$. Indeed, let $\pi_\gamma$ be the invariant probability measure for the
process governed by the operator $L_y + \gamma \mathcal{D}_y$ on $\mathbf{S}$. Integrating the equality 
$ M(\gamma) \varphi_\gamma = \lambda_\gamma \varphi_\gamma $ with respect to $\pi_\gamma$, and dividing both sided by $\gamma$, we obtain
\[
(\gamma - 1) \int_\mathbf{S} \alpha \varphi_\gamma d \pi_\gamma + \int_\mathbf{S} \beta \varphi_\gamma d \pi_\gamma = \gamma^{-1} \lambda_\gamma \int_\mathbf{S} \varphi_\gamma 
d \pi_\gamma.
\]
Since $\alpha > 0$ and $\beta$ is bounded, this implies that $\lim_{\gamma \rightarrow \infty} \lambda_\gamma = +\infty$. Therefore, there exists $\gamma > 0$ such that $\lambda_\gamma = 0$. 

Let us 
show that such $\gamma$ is unique.  Assume the contrary, i.e., that $(\gamma_1, \varphi_1)$ and $(\gamma_2, \varphi_2)$ satisfy 
(\ref{mg}) and $0 < \gamma_1 < \gamma_2$.  (The case when $\gamma_1$ and $\gamma_2$ are negative can
be considered similarly.)  Let $\tilde{X}_t 
= (\tilde{Y}_t, \tilde{Z}_t)$ be the
family of diffusion processes on $\mathbf{S} \times (0, \infty)$ with the generator $K$ defined in (\ref{kok}).
(This operator can be obtained from the generator of  ${X}_t$ by discarding the last term in~(\ref{oprww}).) 

Consider the processes $\xi^1_t =  \varphi_1(\tilde{Y}_t)(\tilde{Z}_t)^{\gamma_1}$ and 
$\xi^2_t = \varphi_2(\tilde{Y}_t)(\tilde{Z}_t)^{\gamma_2}$. Here, we fix an initial point $x$ such that $\xi^1_0 = 1$. Let 
\[
\tau_n = \inf\{t: \xi^1_t = \frac{1}{n}~{\rm or}~ \xi^1_t = n\}.
\]
Since $K (\varphi_1(y) z^{\gamma_1}) = 
K (\varphi_2(y) z^{\gamma_2}) = 0$,  by the Ito formula,
$\xi^1_t$ and $\xi^2_t$ are local martingales, while the stopped processes $\xi^1_{\tau_n \wedge t}$ and $\xi^2_{\tau_n \wedge t}$ are martingales since they are bounded. Note that $\tau_n < \infty$ almost surely since $\gamma_1 
\neq 0$. By the 
Optional Stopping Theorem applied to the process $\xi^1_{\tau_n \wedge t}$, $  \mathrm{E}_x  \xi^1_{\tau_n} = \xi^1_0 = 1$, and
therefore $\mathrm{P}_x(\xi^1_{\tau_n} = 1/n) = n/(n+1)$, while 
$\mathrm{P}_x(\xi^1_{\tau_n} = n) = 1/(n+1)$.  We estimate
\[
\mathrm{E}_x \xi^2_{\tau_n} \geq \inf ({\varphi_2}/{\varphi_1^{\frac{\gamma_2}{\gamma_1}}})  \mathrm{E}_x 
( \xi^1_{\tau_n})^{\frac{\gamma_2}{\gamma_1}} \geq   \inf ({\varphi_2}/{\varphi_1^{\frac{\gamma_2}{\gamma_1}}}) \mathrm{P}_x(\xi^1_{\tau_n} = 
n) n^{\frac{\gamma_2}{\gamma_1}} \rightarrow \infty~~{\rm as}~n \rightarrow \infty,
\]
since $\mathrm{P}_x(\xi^1_{\tau_n} = 
n) n^{{\gamma_2}/{\gamma_1}} = (n+1)^{-1} n^{{\gamma_2}/{\gamma_1}} \rightarrow \infty$.
However, by the 
Optional Stopping Theorem applied to the process $\xi^2_{\tau_n \wedge t}$, $\mathrm{E}_x \xi^2_{\tau_n} = \xi^2_0$ does not depend on $n$.
Thus we get a contradiction, which proves uniqueness. \qed
\medskip
\\ 
{\bf Remark.} If $\beta(y)/\alpha(y) = c$ is constant on $\mathbf{S}$, then equation (\ref{mg}) can be solved explicitly: $\varphi$ is constant and $\gamma = 1 - c$.

\subsection{Transition probabilities near an invariant surface} \label{tpni}

For a closed set $A$, let $\tau(A) = \inf\{t \geq 0: X_t \in A\}$.
 For $\varkappa \geq 0$, we can define 
\begin{equation} \label{ga1}
\Gamma_\varkappa = \{(y,z): (\varphi(y))^{\frac{1}{\gamma}} z = \varkappa\},
\end{equation}
where $\varphi$ is defined in Lemma~\ref{spec}. For $\varkappa_1< \varkappa_2$, we denote the region
between $\Gamma_{\varkappa_1}$ and $\Gamma_{\varkappa_2}$ by
\begin{equation} \label{stri}
V_{\varkappa_1, \varkappa_2} = \{(y,z): \varkappa_1 \leq (\varphi(y))^{\frac{1}{\gamma}} z  \leq  \varkappa_2\}.
\end{equation}

Let us use Lemma~\ref{spec} to study the transition probabilities for the process $X_t$ near $\mathcal{M}$. 
We start with the case of an attracting surface.
\begin{lemma} \label{pratf} Let $\gamma > 0$. For each $\eta > 0$,  for all sufficiently small $\varkappa_1, \varkappa_2$ (depending on $\eta$)   satisfying
$0 < \varkappa_1 < \varkappa_2$,
\begin{equation} \label{prexf}
  \frac{(1 - \eta)\zeta^\gamma - \varkappa_1^\gamma}{\varkappa_2^\gamma - \varkappa_1^\gamma} \leq \mathrm{P}_x 
( X_{ \tau (\Gamma_{\varkappa_1} \bigcup \Gamma_{\varkappa_2})} \in \Gamma_{\varkappa_2}) \leq 
\frac{(1 + \eta) \zeta^\gamma - \varkappa_1^\gamma}{\varkappa_2^\gamma - \varkappa_1^\gamma},
\end{equation}  
provided that $x \in \Gamma_\zeta$ with $ \varkappa_1 \leq \zeta \leq \varkappa_2$.
\end{lemma}
\proof Let us write the process $X_t$ in the coordinate form as $(Y_t, Z_t)$.
The proof relies on the idea that 
the process  $\varphi(Y_t)(Z_t)^\gamma$ is nearly a martingale in $V_{\varkappa_1, \varkappa_2}$ for small $\varkappa_1, \varkappa_2$ (or, equivalently,
the generator of the process applied to $\varphi(y)z^\gamma$ is nearly zero).
Let us now give a rigorous argument. From the proof of Lemma~\ref{spec} and the smooth dependence of the top eigenvalue/eigenfunction
on the parameter $\gamma$, it follows that there exist   $\hat{\gamma} \in (\gamma, \gamma +1)$ and a positive-valued function 
$\hat{\varphi}  \in C^2({\mathbf{S}})$ satisfying $\int_{\mathbf{S}} \hat{\varphi} d \pi = 1$ such that 
\begin{equation} \label{yyu2f}
L_y \hat{\varphi} +  \hat{\gamma} (\hat{\gamma} -1) \alpha \hat{\varphi} +  \hat{\gamma} \beta \hat{\varphi} + \hat{\gamma} \mathcal{D}_y \hat{\varphi} = \hat{c} \hat{\varphi},
\end{equation}
where $\hat{c} >0$. Let
\[
u (y,z) = \varphi(y) z^\gamma,~~ \hat{u} (y,z)  = \hat{\varphi}(y) z^{\hat{\gamma}}.
\]
By Lemma~\ref{geneNN} and (\ref{mg}), the function $f (y,z) = (u (y,z) - \varkappa_1^\gamma)/(\varkappa_2^\gamma - \varkappa_1^\gamma)$ satisfies
\begin{equation} \label{enffp}
L  f  (y,z) = R  f  (y,z),~~~(y,z) \in V_{\varkappa_1, \varkappa_2},
\end{equation}
\[
f|_{\Gamma_{\varkappa_1}} = 0,~~~f|_{\Gamma_{\varkappa_2}} = 1,
\]
where the operator $R$ is defined in Lemma~\ref{geneNN}.
The function $v (x) = \mathrm{P}_x 
( X_{ \tau (\Gamma_{\varkappa_1} \bigcup \Gamma_{\varkappa_2})} \in \Gamma_{\varkappa_2})$, which is what we are interested in, satisfies the same boundary conditions and
the same equation (\ref{enffp}) but with zero instead of $R  f (y,z)$ in the right hand side. Thus, 
\begin{equation} \label{probnf}
\mathrm{P}_x 
( X_{ \tau (\Gamma_{\varkappa_1} \bigcup \Gamma_{\varkappa_2})} \in \Gamma_{\varkappa_2}) = v (x) =  f (x) - g (x),
\end{equation}
 where $g $ solves
\[
L  g  (y,z) = R    f (y,z),~~~(y,z) \in V_{\varkappa_1, \varkappa_2},
\]
\[
g |_{\Gamma_{\varkappa_1}} = 0,~~~g |_{\Gamma_{\varkappa_2}} = 0.
\]
Observe that there is $C > 0$ such that $|R  f  (y,z)| \leq  C  z^{\gamma+1}/(\varkappa_2^\gamma - \varkappa_1^\gamma)   $, provided that $\varkappa_1$ and $\varkappa_2$ are sufficiently small. By Lemma~\ref{geneNN}, using (\ref{yyu2f}), we can find an arbitrarily small 
$\hat{k} > 0$  (by taking  $\varkappa_1, \varkappa_2$
sufficiently small) such that
\[
 \hat{k}  L  \hat{u} (y,z) \geq  C    z^{\gamma+1},~~~(y,z) \in V_{\varkappa_1, \varkappa_2}.
\] 
By taking $\hat{k}$ sufficiently small,  we can find an arbitrarily small $k >0$ such that
\[
k  \inf_{(y,z) \in {\Gamma_{\varkappa_1}} } u(y,z) \geq \hat{k}  \sup _{(y,z) \in {\Gamma_{\varkappa_1}} } \hat{u}(y,z),~~~
k  \inf_{(y,z) \in {\Gamma_{\varkappa_2}} } u(y,z) \geq \hat{k}  \sup _{(y,z) \in {\Gamma_{\varkappa_2}} } \hat{u}(y,z).
\]
 Thus the function $\tilde{g}  =  ( k u -
\hat{k} \hat{u} )/(\varkappa_2^\gamma - \varkappa_1^\gamma)$ satisfies
\[
L\tilde{g} (y,z) \leq  - C (1-k) z^{\gamma+1}/(\varkappa_2^\gamma - \varkappa_1^\gamma),~~~(y,z) \in V_{\varkappa_1, \varkappa_2},
\]
\[
\tilde{g}|_{\Gamma_{\varkappa_1}} \geq 0,~~~\tilde{g}|_{\Gamma_{\varkappa_2}} \geq 0,
\]
which implies that
\[
\tilde{g}(y,z) \geq (1-k) |{g} (y,z)|, ~~~(y,z) \in V_{\varkappa_1, \varkappa_2}.
\]
Thus
\[
 |{g} (y,z)| \leq \frac{ k \varphi(y) z^{\gamma} - \hat{k} \hat{\varphi}(y) z^{\hat{\gamma}}}{(1 - k)(\varkappa_2^\gamma - \varkappa_1^\gamma)}  
\leq \frac{ k \varphi(y) z^{\gamma}}{(1 - k)(\varkappa_2^\gamma - \varkappa_1^\gamma)} .
\]
We are interested in $(y,z) \in \Gamma_{\zeta}$. For each ${\eta} > 0$, by making  $k$ sufficiently small, we can make sure that
\[
\sup_{(y,z) \in \Gamma_{\zeta}} |{g}(y,z)| \leq   \frac{{\eta} \zeta^\gamma}{\varkappa_2^\gamma - \varkappa_1^\gamma}. 
\]
Combined with (\ref{probnf}) and the definition of $f $,  this gives the desired estimate (\ref{prexf}). \qed
\\

We have the following counterpart of Lemma~\ref{pratf} in the case of a repelling surface. It's proof is similar to that of Lemma~\ref{pratf}.
\begin{lemma} \label{prat2} Let $\gamma < 0$. For each $\eta > 0$,  for all sufficiently small $\varkappa_1, \varkappa_2$ 
(depending on $\eta$) satisfying
$0 < \varkappa_1 < \varkappa_2$,
\begin{equation} \label{prex}
  \frac{(1 + \eta)\zeta^\gamma - \varkappa_1^\gamma}{\varkappa_2^\gamma - \varkappa_1^\gamma} \leq \mathrm{P}_x 
( X_{ \tau (\Gamma_{\varkappa_1} \bigcup \Gamma_{\varkappa_2})} \in \Gamma_{\varkappa_2}) \leq 
\frac{(1 - \eta) \zeta^\gamma - \varkappa_1^\gamma}{\varkappa_2^\gamma - \varkappa_1^\gamma},
\end{equation}  
provided that $x \in \Gamma_\zeta$ with $ \varkappa_1 \leq \zeta \leq \varkappa_2$.
\end{lemma}
The following lemma justifies our use of the terms ``attracting'' and ``repelling'' surface. Define
\[
q(x) = \mathrm{P}_x (\lim_{t \rightarrow \infty} {\rm dist}(X_t, \mathcal{M}) = 0). 
\] 
\begin{lemma} \label{cohhi} (a)
If $\gamma > 0$, then $q(x) \rightarrow 1$ as ${\rm dist}(x, \mathcal{M}) \downarrow 0$. 
 (b) If $\gamma < 0$, then
$q(x)= 0$  for each $x \notin \mathcal{M}$.
\end{lemma}
\proof Consider first the case when $\gamma > 0$. Given an arbitrary $\delta > 0$, take $a \in (0,1)$ such that $2 a^\gamma \leq \delta$.
From the second inequality in (\ref{prexf}) with $\eta = 1$ it follows that
\[
 \mathrm{P}_x 
( X_{ \tau (\Gamma_{\varkappa_1} \bigcup \Gamma_{\varkappa_2})} \in \Gamma_{\varkappa_2}) \leq 
\frac{2 (a \varkappa_2)^\gamma - \varkappa_1^\gamma}{\varkappa_2^\gamma - \varkappa_1^\gamma} \leq 2 a^\gamma \leq \delta,
\]
provided that $x \in \Gamma_\zeta$ with $ \varkappa_1 \leq \zeta \leq  a \varkappa_2$. Since $\varkappa_1$ can be taken arbitrarily small, from here it follows that
\[
 \mathrm{P}_x (\liminf_{t \rightarrow \infty} {\rm dist}(X_t, \mathcal{M}) = 0) \geq 1 - \delta.
\]
Again, from the second inequality in (\ref{prexf}) it follows that $\liminf_{t \rightarrow \infty} {\rm dist}(X_t, \mathcal{M}) = 0$ implies (except on an event of probability zero) that $\lim_{t \rightarrow \infty} {\rm dist}(X_t, \mathcal{M}) = 0$. Hence, 
\[
 \mathrm{P}_x (\lim_{t \rightarrow \infty} {\rm dist}(X_t, \mathcal{M}) = 0) \geq 1 - \delta
\]
for $x \in \Gamma_\zeta$ with sufficiently small $  \zeta$ (depending on $\delta$). 
Since $\delta$ was arbitrary, the first statement of the lemma follows. 

In order to prove the second statement, take $\varkappa_2 >0$ sufficiently small so that the first inequality in (\ref{prex}) is valid with 
$\eta = 1$. For each $x$ such that $x \in \Gamma_\zeta$ with $\zeta \leq \varkappa_2$ and each $\varkappa_1 < \zeta$,
\[
  \frac{2\zeta^\gamma - \varkappa_1^\gamma}{\varkappa_2^\gamma - \varkappa_1^\gamma} \leq \mathrm{P}_x 
( X_{ \tau (\Gamma_{\varkappa_1} \bigcup \Gamma_{\varkappa_2})} \in \Gamma_{\varkappa_2}).
\]
Since $\gamma < 0$ now,  taking $\varkappa_1 \downarrow 0$, we obtain  $\mathrm{P}_x 
( \tau (\Gamma_{\varkappa_2}) < \infty) = 1$. Since this holds for each $x$ such that $x \in \Gamma_\zeta$ with $\zeta \leq \varkappa_2$, 
we conclude that $\mathrm{P}_x (\lim_{t \rightarrow \infty} {\rm dist}(X_t, \mathcal{M}) = 0) = 0$. \qed
\\

The next lemma describes the distribution of $Y_t$ component of the process, provided that the process starts close to $\mathcal{M}$.
Recall that $\pi_{\mathcal{M}}$ is the invariant measure for the process on $\mathcal{M}$.  
\begin{lemma} \label{proxnnn}
For each $f \in C_b(\mathbb{R}^d)$ and   each $\delta > 0$,  there exist $\varkappa > 0$ and $t > 0$ such that
\[
|\mathrm{E}_x f(X_t) - \int_{\mathcal{M}} {f} d\pi_{\mathcal{M}}| < \delta
\]
when $x \in V_{0,\varkappa}$.
\end{lemma}
\proof 
Let $p: \mathcal{M}^r \rightarrow \mathcal{M}$ be the projection of an $r$-neighborhood of $\mathcal{M}$ (with sufficiently small 
$r$)   onto $\mathcal{M}$, 
i.e., $p(x) = m$ for each $x = (m, n, z) \in \mathcal{M}^r$.
 By Assumption (a), we can find $t > 0$ such that 
\[
|\mathrm{E}_x f(X_t) - \int_{\mathcal{M}} {f} d\pi_{\mathcal{M}}| <  \delta/2.
\]
for each $x \in \mathcal{M}$. Since $f$ is bounded, continuous, and $\|x - p(x)\| \rightarrow 0$ as 
${\rm dist}(x, \mathcal{M}) \rightarrow 0$, we can find
 $\varkappa > 0$ such that  
\begin{equation} \label{pron1}
|\mathrm{E}_x f(X_{t  }) - \mathrm{E}_{p(x)} f(X_{t  })| < \delta/2,~~~x \in V_{0, \varkappa}.
\end{equation}
Together with the previous inequality, this implies the  desired result.
\qed

\subsection{Limiting distribution of the unperturbed process} \label{ldup}
In the previous sections, we considered the behavior of the process near a single invariant surface. Recall that, in general, there 
may be multiple such surfaces $\mathcal{M}_1,...,\mathcal{M}_m$. Let us assume that $d - d_k > 1$ for each $k =1,...,m$. The situation when
$d_k = d-1$ for some $k$ can be considered similarly by treating each connected component of $D$ separately. Recall that, 
by Lemma~\ref{spec}, there is a constant $\gamma_k$ associated to each surface $\mathcal{M}_k$. Without loss of generality, we assume that
$\gamma_1 \geq \gamma_2 \geq ... \geq \gamma_m$. 

Consider first the case when at least one surface is attracting, i.e., 
$\gamma_1 > 0$. Let $\overline{m}$ be such that $\gamma_1 \geq ... \geq \gamma_{\overline{m}} > 0 > \gamma_{\overline{m}+1}
\geq  ... \geq \gamma_m$. For  $k \leq \overline{m}$, let $E_k$ be the event that $\lim_{t \rightarrow \infty}  
{\rm dist}(X_t, \mathcal{M}_k) = 0$.  Define $p^x_k = \mathrm{P}_x(E_k)$. 
Let $\mathbf{\pi}_{\mathcal{M}_k}$ be the marginal distribution of $\pi_k$ on $\mathcal{M}_k$ (see the discussion after
Lemma~\ref{geneNN}).  
\begin{theorem} \label{unpnn} Let $\gamma_1 > 0$.
The functions $p^x_k$ are positive and continuous in $x \in D \bigcup \mathcal{M}_1 \bigcup... \bigcup \mathcal{M}_{\overline{m}}$ for each $1 \leq k \leq \overline{m}$ and satisfy $\sum_{k=1}^{\overline{m}} p^x_k = 1$.
The distribution of the process $X^x_t$ converges weakly, as $t \rightarrow \infty$, to the measure $\sum_{k=1}^{\overline{m}} 
p^x_k \pi_{\mathcal{M}_k}$. 
\end{theorem}
\proof  Define $\Gamma^k_\varkappa$ and $V^k_{\varkappa_1, \varkappa_2}$, $1 \leq k \leq m$,  as in (\ref{ga1}) and (\ref{stri}) but
for the surface $\mathcal{M}_k$ instead of a generic surface $\mathcal{M}$. For $\varkappa > 0$, let $E_{k,\varkappa}$, $1 \leq k 
\leq \overline{m}$,
 be the event that $\tau(V^k_{0,\varkappa}) = \min(\tau(V^1_{0,\varkappa}),...,  \tau(V^{\overline{m}}_{0,\varkappa}))$.
 Define $p^x_{k,\varkappa} = \mathrm{P}_x(E_{k,\varkappa})$. Since the process is non-degenerate in $D$, does not escape to infinity (assumption (d)),  and leaves any sufficiently small 
neighborhood of $\mathcal{M}_{\overline{m}+1}\bigcup ... \bigcup \mathcal{M}_m$ with probability one (Lemma~\ref{cohhi}, part (b)), 
the functions $p^x_{k,\varkappa}$, $1 \leq k \leq \overline{m}$, are continuous on $D \bigcup \mathcal{M}_1 \bigcup... \bigcup \mathcal{M}_{\overline{m}}$,  and satisfy $\sum_{k=1}^{\overline{m}} p^x_{k,
\varkappa} = 1$ there. 

By part (a) of Lemma~\ref{cohhi} and the strong Markov property, $p^x_{k, \varkappa} \rightarrow p^x_k$ as $\varkappa \downarrow 0$. The convergence is uniform on any compact subset of $D \bigcup \mathcal{M}_1 \bigcup... \bigcup \mathcal{M}_{\overline{m}}$. This implies that $p^x_k$ are continuous in $x \in D$ for each $1 \leq k \leq \overline{m}$ and satisfy $\sum_{k=1}^{\overline{m}} p^x_k = 1$. The positivity of $p^x_k$ follows from the non-degeneracy of the process $X_t$ in
$D$. The convergence of the distribution of $X^x_t$   
to the measure $\sum_{k=1}^{\overline{m}}  p^x_k \pi_{\mathcal{M}_k}$ now follows from Lemma~\ref{proxnnn} and the Markov property of the
process. 
\qed
\\

Next, we consider the case when all the  components of the boundary are repelling. 
\begin{theorem} \label{repth} Let $\gamma_1 < 0$.  Suppose that $d - d_k > 1$ for each $k =1,...,m$, i.e., $D$ is connected.  There is a unique measure $\mu$ on $D$ that is invariant for the
 process $X_t$. This measure
is absolutely continuous with respect to the Lebesgue measure. The
the distribution of $X_{t }$ converges weakly, as $t \rightarrow \infty$, to $\mu$ uniformly with respect to initial point on a given compact set. Namely, for each compact $K \subset D$ and each $f \in C_b(D)$, 
\begin{equation} \label{wko}
\lim_{t \rightarrow \infty} \sup_{x \in K} |\mathrm{E}_x f(X_{t}) - \int_D f d\mu| =0.
\end{equation}
\end{theorem}
\proof We will construct a function $V \in C^2(D)$ such that (a) $V(x) >0$, $x \in D$, (b)~$\lim_{x \rightarrow \mathcal{M}_k} V(x) 
= \infty$, $1 \leq k \leq m$,  $\lim_{\|x\| \rightarrow \infty} V(x) = \infty$, (c) $L V \leq K - c V$ for some positive constants $K$, $c$. 
Since the process is non-degenerate in $D$, by \cite{Har}   there exists a probability measure $\mu$ on $D$ that is invariant with respect to
the Lebesgue measure and satisfies, for some $C > 0$ and $0 < \rho < 1$, 
\[
\sup_{x \in D}\frac{|\mathrm{E}_x f(X_t) -\int_D f d \mu | }{V(x)} \leq C \rho^t   
\sup_{x \in D}\frac{|f(x) -\int_D f d \mu | }{V(x)},~~t \geq 0,
\] 
whenever $f \in C(D)$ is such that the right-hand side is finite. This easily implies (\ref{wko}).

It remains to exhibit the function $V$. Consider a neighborhood a surface  $\mathcal{M} = \mathcal{M}_k$ and temporarily 
drop the subscript $k$ from the notation. From the proof of Lemma~\ref{spec} and the smooth dependence of the top eigenvalue/eigenfunction
on the parameter $\gamma$, it follows that there exist   $\hat{\gamma} \in (\gamma, 0)$ and a positive-valued function 
$\hat{\varphi}  \in C^2({\mathbf{S}})$ satisfying $\int_{\mathbf{S}} \hat{\varphi} d \pi = 1$ such that 
\[
L_y \hat{\varphi} +  \hat{\gamma} (\hat{\gamma} -1)\alpha  \hat{\varphi} +  \hat{\gamma} \beta \hat{\varphi} + \hat{\gamma} \mathcal{D}_y \hat{\varphi} = \hat{c} \hat{\varphi},
\]
where $\hat{c} < 0$. Define $V (y,z)  = \hat{\varphi}(y) z^{\hat{\gamma}}$ in a small neighborhood of $\mathcal{M}$. By Lemma~\ref{geneNN},
\[
LV (y,z) = \hat{c} V(y,z) + R V(y,z)  \leq  \frac{\hat{c}}{2} V(y,z),
\]
where $R$ is the same as in  Lemma~\ref{geneNN} and the inequality holds if $z$ is sufficiently small. In this way, $V$ is defined in a neighborhood of each of the surfaces $\mathcal{M}_k$. For $\|x\|$ sufficiently large, we define $V(x) = \exp(\delta \|x\|)$. It follows from
assumption (d) on the vector fields (see Section~\ref{intro}) that $L V(x) \leq - c V(x)$ if $\|x\|$ is sufficiently large, $\delta > 0$ is sufficiently small, and $c > 0$ is sufficiently small. Having defined $V$ in the neighborhoods of $\mathcal{M}_k$ and a neighborhood of 
infinity, we can extend it to $D$ as a positive $C^2$ function.  Thus $L V \leq K - c V$ for some positive constants $K$, $c$. \qed
\\

\section{Asymptotic behavior of the perturbed process}
We will now examine the asymptotic behavior of the process $X^{\varepsilon}_t$ defined in Section~\ref{intro}. While some of the analysis 
from the previous section is still applicable, there are a couple of new technical issues. First, estimates on the transition probabilities are now valid only away from $\varepsilon$-dependent neighborhoods of the invariant surfaces due to the presence of a small
non-degenerate perturbation. Second, the perturbed process will eventually escape from a neighborhood of an invariant surface even if the 
surface is attracting. Estimates on the typical escape times are of crucial importance for the description of the metastable behavior.  

\subsection{Transition probabilities  near an invariant surface} \label{tpni2}
We start with the counterpart to Lemma~\ref{pratf} for the perturbed process. For a closed set $A$, 
we define $\tau^{\varepsilon}(A) = \inf\{t \geq 0: X^{\varepsilon}_t \in A\}$.
\begin{lemma} \label{pratfzzz} Let $\gamma > 0$. For each $\eta > 0$,  for all sufficiently large $r$ (depending on $\eta$), for all sufficiently small $\varkappa_2 >0$ (depending on $\eta$), for all sufficiently small $\varepsilon$ 
(depending on $r$ and $\varkappa_2$), and for all $\varkappa_1 \in [r \varepsilon, \varkappa_2)$,  
\begin{equation} \label{prexfzxz}
  \frac{(1 - \eta)\zeta^\gamma - \varkappa_1^\gamma}{\varkappa_2^\gamma - \varkappa_1^\gamma} \leq \mathrm{P}_x 
( X^{\varepsilon}_{ \tau^{\varepsilon} (\Gamma_{\varkappa_1} \bigcup \Gamma_{\varkappa_2})} \in \Gamma_{\varkappa_2}) \leq 
\frac{(1 + \eta) \zeta^\gamma - \varkappa_1^\gamma}{\varkappa_2^\gamma - \varkappa_1^\gamma},
\end{equation}  
provided that $x \in \Gamma_\zeta$ with $ \varkappa_1 \leq \zeta \leq \varkappa_2$.
\end{lemma}
\proof The proof is similar to that of Lemma~\ref{pratf}. However, an extra difficulty appears due to the presence of the term 
$ \varepsilon^2   \tilde{L} f$ in the right-hand side of the analogue of  (\ref{enffp}) if $L$ is replaced by $L^\varepsilon$. The following observation will help us cope with the extra term: for each $\gamma' > 0$,  it follows by re-writing $\tilde{L}$ in $(y,z)$-coordinates that 
\begin{equation} \label{zzui}
| \tilde{L} z^{\gamma'}| \leq K   |z^{\gamma'-2}|,
\end{equation}
for some constant $K$, provided that $|z| > 0$ is sufficiently small.  

Recall from the proof of Lemma~\ref{spec} that the top eigenvalue of the operator 
in the left hand side of (\ref{mg}) satisfies $\lambda_0 = 0$, $\lambda'_0 < 0$, and $\lambda_\gamma = 0$. Since the top
eigenvalue depends continuously on the parameter,  by the uniqueness part of Lemma~\ref{spec}, the eigenvalues corresponding to the values of the parameter that are slightly 
smaller than $\gamma$ are negative, and thus
 there exist $\gamma_1 \in (\max(0, \gamma-1) ,\gamma)$ and a positive-valued function 
$\varphi_1  \in C^2(S)$ satisfying $\int_S \varphi_1 d \pi = 1$ such that
\begin{equation} \label{yyu1}
L_y \varphi_1 +  \gamma_1 (\gamma_1 -1) \alpha \varphi_1 +  \gamma_1 \beta \varphi_1 + \gamma_1 \mathcal{D}_y \varphi_1 = -c_1 \varphi_1,
\end{equation}
where $c_1 >0$.  Similarly, there exist  $\gamma_2 \in (\gamma, \gamma +1)$ and a positive-valued function 
$\varphi_2  \in C^2(S)$ satisfying $\int_S \varphi_2 d \pi = 1$ such that
\begin{equation} \label{yyu2}
L_y \varphi_2 +  \gamma_2 (\gamma_2 -1) \alpha \varphi_2 +  \gamma_2 \beta \varphi_2 + \gamma_2 \mathcal{D}_y \varphi_2 = c_2 \varphi_2,
\end{equation}
where $c_2 >0$. Let
\[
u(y,z) = \varphi(y) z^\gamma,~~ u_1(y,z)  = \varphi_1(y) z^{\gamma_1},~~u_2(y,z)  = \varphi_2(y) z^{\gamma_2}.
\]
By Lemma~\ref{geneNN} and (\ref{mg}), the function 
\begin{equation} \label{fnff}
f (y,z) = (u(y,z) - \varkappa_1^\gamma)/(\varkappa_2^\gamma - \varkappa_1^\gamma)
\end{equation}
 satisfies
\[
L^\varepsilon f  (y,z) = (R +\varepsilon^2   \tilde{L})f  (y,z),~~~(y,z) \in V_{\varkappa_1, \varkappa_2},
\]
\[
f|_{\Gamma_{\varkappa_1}} = 0,~~~f|_{\Gamma_{\varkappa_2}} = 1.
\]
The function $v^\varepsilon(x) = \mathrm{P}_x 
( X^{\varepsilon}_{ \tau^{\varepsilon} (\Gamma_{\varkappa_1} \bigcup \Gamma_{\varkappa_2})} \in \Gamma_{\varkappa_2})$, which is what we are interested in, satisfies the same boundary conditions and
the same equation but with zero instead of $(R +\varepsilon^2   \tilde{L})f (y,z)$ in the right hand side. Thus, 
\begin{equation} \label{probn}
\mathrm{P}_x 
( X^{\varepsilon}_{ \tau^{\varepsilon} (\Gamma_{\varkappa_1} \bigcup \Gamma_{\varkappa_2})} \in \Gamma_{\varkappa_2}) = 
v^\varepsilon(x) =  
f (x) - g^\varepsilon(x),
\end{equation}
 where $g^\varepsilon$ solves
\begin{equation} \label{sl11}
L^\varepsilon g^\varepsilon (y,z) = (R +\varepsilon^2   \tilde{L})f  (y,z),~~~(y,z) \in V_{\varkappa_1, \varkappa_2},
\end{equation}
\[
g^\varepsilon|_{\Gamma_{\varkappa_1}} = 0,~~~g^\varepsilon|_{\Gamma_{\varkappa_2}} = 0.
\]
Observe that, by (\ref{zzui}), there is $C > 0$ such that 
\begin{equation} \label{rs1s}
|(R +\varepsilon^2   \tilde{L})f  (y,z)| \leq  C   (z^{\gamma+1} + \varepsilon^2 z^{\gamma-2})/(\varkappa_2^\gamma - \varkappa_1^\gamma)
\end{equation}
for all sufficiently small $\varepsilon$. By Lemma~\ref{geneNN}, using (\ref{yyu1}) and (\ref{zzui}), we can find an arbitrarily small $k_1 > 0$  such that,
for all  sufficiently large $r$ (depending on $k_1$), all sufficiently small $\varkappa_2$ (depending on $k_1$), 
\[
k_1 \varepsilon^{ \gamma  - \gamma_1 } L^\varepsilon ( u_1(y,z)) \leq - C   \varepsilon^2 z^{ \gamma-2},~~~(y,z) \in V_{r \varepsilon, \varkappa_2},
\] 
holds for all sufficiently small $\varepsilon$. Similarly,  using (\ref{yyu2}) and (\ref{zzui}), we can find an arbitrarily small $k_2 > 0$  such that,
for all  sufficiently large $r$ (depending on $k_2$), all sufficiently small $\varkappa_2$ (depending on $k_2$), 
\[
 k_2  L^\varepsilon ( u_2(y,z)) \geq  C    z^{\gamma+1},~~~(y,z) \in V_{r \varepsilon, \varkappa_2},
\] 
holds for all sufficiently small $\varepsilon$. By taking $k_2$ sufficiently small,  we can find an arbitrarily small $k >0$ such that
\[
k  \inf_{(y,z) \in {\Gamma_{\varkappa_1}} } u(y,z) \geq k_2  \sup _{(y,z) \in {\Gamma_{\varkappa_1}} } u_2(y,z),~~~
k  \inf_{(y,z) \in {\Gamma_{\varkappa_2}} } u(y,z) \geq k_2  \sup _{(y,z) \in {\Gamma_{\varkappa_2}} } u_2(y,z).
\]
 Thus, for $\varkappa_1 \in [r \varepsilon, \varkappa_2)$, the function $\tilde{g}^\varepsilon  = (k_1 \varepsilon^{
\gamma  - \gamma_1} u_1 + k u -
k_2 u_2)/(\varkappa_2^\gamma - \varkappa_1^\gamma)$ satisfies
\[
L^\varepsilon \tilde{g}^\varepsilon (y,z) \leq  - C (1-k)  (z^{\gamma+1} + \varepsilon^2 z^{ \gamma-2 })/(\varkappa_2^\gamma - \varkappa_1^\gamma),~~~(y,z) \in V_{\varkappa_1, \varkappa_2},
\]
\[
\tilde{g}^\varepsilon|_{\Gamma_{\varkappa_1}} \geq 0,~~~\tilde{g}^\varepsilon|_{\Gamma_{\varkappa_2}} \geq 0.
\]
Comparing this with (\ref{sl11}), (\ref{rs1s}) and using the stochastic representation for the solutions $\tilde{g}^\varepsilon$ and ${g}^\varepsilon$
of the respective equations, we obtain that 
\[
\tilde{g}^\varepsilon(y,z) \geq (1-k) |{g}^\varepsilon(y,z)|, ~~~(y,z) \in V_{\varkappa_1, \varkappa_2}.
\]
Thus
\[
 |{g}^\varepsilon(y,z)| \leq \frac{k_1 \varepsilon^{\gamma  - \gamma_1} \varphi_1(y) z^{\gamma_1} + k \varphi(y) z^{\gamma} - k_2 \varphi_2(y) z^{\gamma_2}}{(1 - k)(\varkappa_2^\gamma - \varkappa_1^\gamma)}  
\leq  \frac{k_1 \varepsilon^{\gamma  - \gamma_1} \varphi_1(y) z^{\gamma_1} + k \varphi(y) z^{\gamma} }{(1 - k)(\varkappa_2^\gamma - \varkappa_1^\gamma)}  .
\]
For each ${\eta} > 0$, by making $k_1$ and $k$ sufficiently small, we can make sure that
\[
\sup_{(y,z) \in \Gamma_{\zeta}} |{g}^\varepsilon(y,z)| \leq   \frac{{\eta} \zeta^\gamma}{\varkappa_2^\gamma - \varkappa_1^\gamma}. 
\]
 Combining this with (\ref{probn}) and the definition (\ref{fnff}) of $f $,
we get the desired result.  \qed
\\

We state the counterpart of Lemma~\ref{pratfzzz} in the case of a repelling surface. It can be
proved similarly to Lemma~\ref{pratfzzz}.
\begin{lemma} \label{prat2rt} Let  $\gamma < 0$. For each $\eta > 0$,  for all sufficiently large $r$ (depending on $\eta$), for all sufficiently small $\varkappa_2 >0$ (depending on $\eta$), for all sufficiently small $\varepsilon$ (depending on $r$ and $\varkappa_2$), and for all $\varkappa_1 \in [r \varepsilon, \varkappa_2)$,
\[
  \frac{(1 + \eta)\zeta^\gamma - \varkappa_1^\gamma}{\varkappa_2^\gamma - \varkappa_1^\gamma} \leq \mathrm{P}_x 
( X^\varepsilon_{ \tau (\Gamma_{\varkappa_1} \bigcup \Gamma_{\varkappa_2})} \in \Gamma_{\varkappa_2}) \leq 
\frac{(1 - \eta) \zeta^\gamma - \varkappa_1^\gamma}{\varkappa_2^\gamma - \varkappa_1^\gamma},
\] 
provided that $x \in \Gamma_\zeta$ with $ \varkappa_1 \leq \zeta \leq \varkappa_2$.
\end{lemma}

We also have the following analogue of Lemma~\ref{proxnnn} for the process $X^\varepsilon_t$.
\begin{lemma} \label{proxnnn2}
For each $f \in C_b(\mathbb{R}^d)$ and   each $\delta > 0$,  there exist $\varkappa > 0$, $\varepsilon_0 > 0$, and $t > 0$ such that
\[
|\mathrm{E}_x f(X^\varepsilon_t) - \int_{\mathcal{M}} {f} d\pi_{\mathcal{M}}| < \delta
\]
when $x \in V_{0,\varkappa}$ and $0 < \varepsilon \leq \varepsilon_0$.
\end{lemma}
The proof of this lemma is similar to that of Lemma~\ref{proxnnn}. The only difference is that, instead of (\ref{pron1}), we can use the
fact that
there exist $\varkappa > 0$ and $\varepsilon_0$ such that  
\[
|\mathrm{E}_x f(X^\varepsilon_{t  }) - \mathrm{E}_{p(x)} f(X_{t  })| < \delta/2,~~~x \in V_{0, \varkappa},~~0 < \varepsilon \leq \varepsilon_0.
\] 
 
%
%
%

\subsection{Transition times near an invariant surface} \label{sttim}
The following lemma bounds the expectation of the time required to exit an $\varepsilon$-dependent neighborhood of an attracting or repelling surface. 
\begin{lemma} \label{epsde}
For each $r >0$, there are $a > 0$ and $\varepsilon_0 > 0$ (that depend on $r$) such that
\[
\mathrm{E}_x \tau^{\varepsilon} (\Gamma_{r \varepsilon}) \leq a
\]
for all $x \in V_{0, r\varepsilon}$, $0 < \varepsilon \leq \varepsilon_0$.
\end{lemma}
\proof Recall that the
generator of the process $X^{\varepsilon}_t$ is the operator $L^\varepsilon = L + \varepsilon^2 \tilde{L}$, where
$L$ can be expressed in $(y,z)$ coordinates using (\ref{oprww}), while $\tilde{L}$ is elliptic. Let $r_0 = {r}/{\inf_{y \in \mathbf{S}} \varphi(y)^{\frac{1}{\gamma}}}$.
Thus $z \leq r_0 \varepsilon $ if $(y,z) \in V_{0, r\varepsilon}$ (see (\ref{stri})).
Define the functions 
\[
f(z) = \exp( c z^2),~~\Phi(z) = \int_0^z f^{-1} (s)ds,~~F(z) = -\int_0^z \Phi(s)f(s) ds,~~ 
\]
\[
 h_\varepsilon(z) =F(z/\varepsilon) - F(r_0),
\] 
where $c > 0$ is a parameter to be chosen later. Note that $h_\varepsilon$ can be viewed as a function on $V_{0, r\varepsilon}$ that is
independent of the $y$ variable.
From (\ref{oprww}) it follows that, for all sufficiently small $\varepsilon$, 
\begin{equation} \label{ft1e}
 L h_\varepsilon (y,z)  \leq 2 \varepsilon^{-1}   (\sup_{y \in \mathbf{S}} |\beta(y)|+1) z\Phi(z) f(z),~~(y,z) \in  V_{0, r\varepsilon},
\end{equation}
where $\beta$ is the coefficient at the first order derivative in (\ref{oprww}). Let us write the operator $\tilde{L}$ in $(y,z)$ coordinates
as 
\[
\tilde{L} u(y,z) = a(y,z) \frac{\partial^2 u(y,z)}{\partial z^2} + b(y,z) \frac{\partial u(y,z)}{\partial z} + \tilde{K} u(y,z),
\]
where $a(y,z) \geq \overline{a} > 0$ and $|b(y,z)| \leq \overline{b}$, $(y,z) \in V_{0, r\varepsilon}$, $a$ and $b$ are some continuous functions, and $\tilde{K}$ is an operator that
includes partial derivatives in $y$ and mixed derivatives and therefore satisfies $\tilde{K}h_\varepsilon =0$. Thus
\begin{equation} \label{ft2e}
\varepsilon^2 \tilde{L} h_\varepsilon (y,z) \leq - 2 \overline{a} c  \varepsilon^{-1}    z\Phi(z) f(z) - \overline{a} + \varepsilon \overline{b} |\Phi(z)|
f(z),~~(y,z) \in  V_{0, r\varepsilon}.
\end{equation}
Combining (\ref{ft1e}) with (\ref{ft2e}) and selecting a sufficiently large $c$ and a sufficiently small $\varepsilon_0$, we obtain
\[
L^\varepsilon h_\varepsilon (y,z) = (L + \varepsilon^2 \tilde{L})h_\varepsilon (y,z) \leq -\overline{a},~~(y,z) \in  V_{0, r\varepsilon},
\]
provided that $0 < \varepsilon \leq \varepsilon_0$. Observe also that 
\[
h_\varepsilon (z) \geq 0,~~{\rm whenever}~~(y,z) \in \partial V_{0, r\varepsilon} = \Gamma_{r \varepsilon}.
\]
On the other hand, the function $\tilde{h}_\varepsilon(y,z) := \mathrm{E}_x \tau^{\varepsilon} (\Gamma_{r \varepsilon}) $ satisfies 
\[
L^\varepsilon \tilde{h}_\varepsilon  (y,z) = -1,~~(y,z) \in  V_{0, r\varepsilon};~~~~ \tilde{h}_\varepsilon|_{\Gamma_{r \varepsilon}} \equiv 0.
\]
From here we conclude that $\tilde{h}_\varepsilon(y,z) \leq h_\varepsilon(z)/\overline{a}$ for $(y,z) \in V_{0, r\varepsilon}$.
This implies the statement of the lemma with $a = -F(r_0) \overline{a}$. 
\qed
\\
\begin{lemma} \label{logexit1b}
(a)  Let $\gamma > 0$. There is $a > 0$ such that for all sufficiently large $r > 0$ and all 
sufficiently small $\varkappa_2$ there is $\varepsilon_0 > 0$ (that depends on $r$ and $\varkappa_2$) such that for all  $\varkappa_1$ satisfying
$r \varepsilon \leq \varkappa_1 < \varkappa_2$,  we have
\[
\mathrm{E}_x \tau^{\varepsilon}(\Gamma_{\varkappa_1} \bigcup \Gamma_{\varkappa_2}) \leq a (1 + \ln(\zeta(x)/\varkappa_1)),~~~x \in V_{\varkappa_1,\varkappa_2},~~0 < \varepsilon \leq \varepsilon_0,
\]
where $\zeta = \zeta(x)$ is such that $x \in \Gamma_\zeta$. 

(b) Let $\gamma < 0$.  There is $b >0$ and such that for all sufficiently small $\varkappa$ there is $\varepsilon_0 > 0$ (that depends on $\varkappa$) such that 
\begin{equation} \label{firmo}
\mathrm{E}_x \tau^{\varepsilon}( \Gamma_{\varkappa}) \leq b (1 + {\rm min} (\ln(\varkappa/\varepsilon), \ln(\varkappa/\zeta) ),~~~x \in V_{0,\varkappa},~~
0 < \varepsilon \leq \varepsilon_0,
\end{equation}
where $\zeta = \zeta(x)$ is such that $x \in \Gamma_\zeta$.
\end{lemma}
\proof Consider the case when $\gamma > 0$ first. As in the proof of Lemma~\ref{pratfzzz}, there exist $\gamma_1 \in (\max(0, \gamma-1) ,\gamma)$ and a positive-valued function  $\varphi_1  \in C^2(S)$ satisfying $\int_S \varphi_1 d \pi = 1$  and (\ref{yyu1}). Let
$u_1(y,z)  = \varphi_1(y) z^{\gamma_1}$. 
By Lemma~\ref{geneNN}, using (\ref{yyu1}), 
for all  sufficiently large $r$ and all sufficiently small $\varkappa_2$, 
\[
L^\varepsilon u_1(y,z) \leq - \frac{c_1}{2} u_1(y,z),~~~(y,z) \in V_{r \varepsilon, \varkappa_2},
\] 
holds for all sufficiently small $\varepsilon$. Let $w(y,z)  = \ln(u_1(y,z)/(\varkappa_1^{\gamma_1} \inf \varphi_1))$,
where $r \varepsilon \leq \varkappa_1 < \varkappa_2$.  
Thus
\[
L^\varepsilon w(y,z) \leq - \frac{c_1}{2},~~~(y,z) \in V_{\varkappa_1, \varkappa_2}.
\]
Since $w(x) \geq 0$ for $x \in V_{\varkappa_1, \varkappa_2}$, 
\[
\mathrm{E}_x \tau^{\varepsilon}(\Gamma_{\varkappa_1} \bigcup \Gamma_{\varkappa_2}) \leq \frac{2}{c_1} w(x),~~~
x \in V_{\varkappa_1, \varkappa_2}.
\]
This implies part (a) of the lemma.

Now consider the case when $\gamma < 0$. Similarly to the way it was done in part (a), it is not difficult to show that there is
$a >0$ such that for all sufficiently large  $r > 0$ and sufficiently small $\varkappa_2 >0$  there is $\varepsilon_0 > 0$ 
such that for all $\varkappa_1 $ satisfying
$r \varepsilon \leq \varkappa_1 < \varkappa_2$,  we have
\[
\mathrm{E}_x \tau^{\varepsilon}(\Gamma_{\varkappa_1} \bigcup \Gamma_{\varkappa_2}) \leq a (1 + \ln(\varkappa_2/\zeta(x))),~~~x \in V_{\varkappa_1,\varkappa_2},~~0 < \varepsilon \leq \varepsilon_0,
\]
where $\zeta = \zeta(x)$ is such that $x \in \Gamma_\zeta$. In particular, with $\varkappa_1 = r \varepsilon$,  we obtain
\begin{equation} \label{piioi}
\mathrm{E}_x \tau^{\varepsilon}(\Gamma_{r \varepsilon} \bigcup \Gamma_{\varkappa_2}) \leq a (1 + \ln(\varkappa_2/\zeta(x))),~~~
x \in V_{r \varepsilon,\varkappa_2},~~0 < \varepsilon \leq \varepsilon_0.
\end{equation}
Pick $\eta > 0$ such that $(1 +\eta)2^\gamma < 1$. By Lemma~\ref{prat2rt}, we can choose $r$ sufficiently large and $\varkappa_2$ sufficiently small so that
\begin{equation} \label{pkko}
\mathrm{P}_x 
( X^\varepsilon_{ \tau (\Gamma_{r \varepsilon} \bigcup \Gamma_{\varkappa_2})} \in \Gamma_{\varkappa_2}) \geq   \frac{(1 + \eta)\zeta^\gamma - (r \varepsilon)^\gamma}{\varkappa_2^\gamma - (r \varepsilon)^\gamma},
\end{equation}
provided that $x \in \Gamma_\zeta$ with $ r \varepsilon \leq \zeta \leq \varkappa_2$.
In particular, there is $c > 0$ such that
\[
  \mathrm{P}_x 
( X^\varepsilon_{ \tau (\Gamma_{r \varepsilon} \bigcup \Gamma_{\varkappa_2})} \in \Gamma_{\varkappa_2}) \geq 
c,
\]
provided that $x \in \Gamma_{2 r \varepsilon}$ and $\varepsilon$ is sufficiently small. By the strong Markov property (considering the transition from $V_{0, r \varepsilon}$ to
$\Gamma_{2 r \varepsilon}$ and then  to $\Gamma_{r \varepsilon} \bigcup \Gamma_{\varkappa_2}$), 
\[
\sup_{x \in V_{0, r\varepsilon}} \mathrm{E}_x \tau^{\varepsilon}( \Gamma_{\varkappa_2}) \leq 
\sup_{x \in V_{0, r\varepsilon}} \mathrm{E}_x \tau^{\varepsilon}( \Gamma_{2 r \varepsilon}) + \sup_{x \in V_{r \varepsilon, \varkappa_2}}\mathrm{E}_x \tau^{\varepsilon}(\Gamma_{r \varepsilon} \bigcup \Gamma_{\varkappa_2}) + (1-c) \sup_{x \in V_{0, r\varepsilon}} \mathrm{E}_x \tau^{\varepsilon}( \Gamma_{\varkappa_2}). 
\]

By Lemma~\ref{epsde}, the first term on the right-hand side is bounded by a constant. Therefore, by (\ref{piioi}), there is a constant $a_1 > 0$ such that, for sufficiently large $r$ and sufficiently small $\varkappa_2$,
\[
c \sup_{x \in V_{0, r\varepsilon}} \mathrm{E}_x \tau^{\varepsilon}( \Gamma_{\varkappa_2}) \leq  a_1 (1 + \ln(\varkappa_2/r \varepsilon))
\]
for all sufficiently small $\varepsilon$. Fixing a value of $r$, from here, (\ref{piioi}), (\ref{pkko}), and the strong Markov property, we obtain that there is a constant $a_2$ such that, for all sufficiently small $\varkappa_2$,
\[
 \mathrm{E}_x \tau^{\varepsilon}( \Gamma_{\varkappa_2}) \leq  a_2 (1 + \ln(\varkappa_2/\zeta(x))),
\]
provided that $x \in \Gamma_\zeta$ with $ r \varepsilon \leq \zeta \leq \varkappa_2$. Combining the last two inequalities, we obtain
part (b) of the lemma.
 \qed
\\

A result similar to Lemma~\ref{logexit1b} holds for the unperturbed process. We will need it in the case of a repelling surface. 
We formulate it as a separate lemma. Its proof is completely similar to the proof of part (b) of Lemma~\ref{logexit1b}.
\begin{lemma} \label{ttlrs}
Let $\gamma < 0$.  There is $b >0$ such that, for all sufficiently small $\varkappa > 0$,
\[
\mathrm{E}_x \tau ( \Gamma_{\varkappa}) \leq b (1 + \ln (\varkappa/\zeta(x))),~~~x \in V_{0,\varkappa} \setminus \mathcal{M},
\]
where $\zeta = \zeta(x)$ is such that $x \in \Gamma_\zeta$.
\end{lemma} 

Later, we will use an estimate similar to (\ref{firmo}) but for the second moment of the exit time. In fact, we will only need a crude
version with $x$ restricted to $\Gamma_{\varkappa/2}$. 
\begin{lemma} \label{semom}
Let $\gamma < 0$.  There is $c >0$ and such that for all sufficiently small $\varkappa$ there is $\varepsilon_0 > 0$ (that depends on $\varkappa$) such that 
\[
\mathrm{E}_x (\tau^{\varepsilon}( \Gamma_{\varkappa}))^2 \leq c,~~~x \in \Gamma_{\varkappa/2},~~
0 < \varepsilon \leq \varepsilon_0.
\]
\end{lemma}
\proof Let $u^\varepsilon(x)  = \mathrm{E}_x \tau^{\varepsilon}( \Gamma_{\varkappa})$,   $v^\varepsilon(x) = 
\mathrm{E}_x (\tau^{\varepsilon}( \Gamma_{\varkappa}))^2$, $x \in V_{0,\varkappa}$. Observe that 
\begin{equation} \label{vep}
L^\varepsilon v^\varepsilon (x) = -2 u^\varepsilon (x),~~x \in V_{0,\varkappa},
\end{equation}
 while $v^\varepsilon(x) = 0$, $x \in \Gamma_\varkappa$.

Pick $\hat{\gamma} \in (\gamma, 0)$. Let
$\hat{\varphi} > 0$ and $\hat{c} > 0$ be such that $K (\hat{\varphi}(y) z^{\hat{\gamma}})= - \hat{c} \varphi(y) z^{\hat{\gamma}} $, $y \in \mathbf{S}$, where $K$ is the operator in (\ref{kok}). 
Denote $\hat{v}(x) = \hat{v}(y,z) = \hat{\varphi}(y) z^{\hat{\gamma}}/\varkappa^{\hat{\gamma}}$. Then
\[
L^\varepsilon \hat{v} (x) = K \hat{v} (x) + R \hat{v} (x) + \varepsilon^2 \tilde{L} \hat{v} (x) \leq -\frac{\hat{c}}{2} \hat{v} (x),~~x \in V_{r \varepsilon, \varkappa},
\]
provided that   $r$ is sufficiently large  and $\varkappa, \varepsilon$ are sufficiently small. We can compare the right hand side of this inequality with the right hand side of (\ref{vep}):
\[
\frac{\hat{c}\hat{v} (y,z)}{2}  =\frac{\hat{c}\hat{\varphi}(y) z^{\hat{\gamma}} }{2 \varkappa^{\hat{\gamma}}}  \geq k  2 u^\varepsilon (y,z), 
\]
where the inequality holds,  by (\ref{firmo}), for sufficiently small $k > 0$ for $(y,z) \in  V_{r \varepsilon, \varkappa}$. Next, we compare the values of $\hat{v}$ and $v^\varepsilon$ on the boundary of $ V_{r \varepsilon, \varkappa}$. For $x \in \Gamma_\varkappa$, we have $\hat{v}(x) > 0$, while $v^\varepsilon(x) =0$. For $x \in \Gamma_{r \varepsilon}$,  we estimate $v^\varepsilon$ by solving a version of (\ref{vep}) in which the function $u^\varepsilon$ is replaced by its supremum, which, by (\ref{firmo}), is estimated from above by $K \ln(\varkappa/\varepsilon)$ for sufficiently small $\varepsilon$, where $K$ is a positive constant. Thus, using (\ref{firmo}) again, we obtain, with a different constant $K$,
\[
\sup_{x \in \Gamma_{r \varepsilon}} v^\varepsilon(x) \leq K (\ln(\varkappa/\varepsilon))^2 \leq \inf_{x \in \Gamma_{r \varepsilon}} \hat{v}(x),
\]
where the inequalities hold for all sufficiently small $\varepsilon$. Now we can conclude that there is another constant $K > 0$ such that
\[
v^\varepsilon(x) \leq K\hat{v}(x), ~~x \in V_{r \varepsilon, \varkappa},
\]
for all sufficiently small $\varepsilon$. Applying this to $x \in \Gamma_{\varkappa/2}$ completes the proof of the lemma. \qed 
\medskip

%
%

Now we can prove the result concerning the exit time from a fixed neighborhood of an attracting surface.
\begin{lemma} \label{mltime} Let $\gamma > 0$.  

(a) There is $a > 0$ such that for each $r >0$ and each sufficiently small $\varkappa>0$, there is $\varepsilon_0 > 0$ such that
\begin{equation} \label{time1}
   \mathrm{E}_x \tau^{\varepsilon} (\Gamma_\varkappa) \leq a    ({\varkappa}/{\varepsilon})^\gamma,
\end{equation}
for all $x \in V_{0, r\varepsilon}$, $0 < \varepsilon \leq \varepsilon_0$.


(b) For each $r > 0$ and $\delta > 0$, there are $\varepsilon_0   > 0 $ and $s > 0$ such that
\begin{equation} \label{bfb}
\mathrm{P}_x\left( \varepsilon^\gamma   {\tau^{\varepsilon} (\Gamma_\varkappa)}  < s\right) \leq \delta
\end{equation}
for all $x \in V_{0,r\varepsilon}$, 
$0 < \varepsilon \leq \varepsilon_0$. 
\end{lemma} 
\proof We will prove the results for $x \in \Gamma_{r\varepsilon}$, which will imply the results for $x \in V_{0,r\varepsilon}$ since, by
Lemma~\ref{epsde}, $\mathrm{E}_x \tau^{\varepsilon}
(\Gamma_{r \varepsilon})$ is uniformly bounded in   $x \in V_{0,r\varepsilon}$, 
$0 < \varepsilon \leq \varepsilon_0$. Moreover, without loss of generality, we can consider an arbitrarily large value of $r$. In order to 
fix a convenient value of $r$, we first take $\eta$ sufficiently small so that $(1-\eta)2^\gamma > 1$. Then the lower bound in 
(\ref{prexfzxz}) is nontrivial if $\varkappa_1 = r\varepsilon$ and $\zeta = 2 r \varepsilon$. We then fix a sufficiently large value of $r$ so
that (\ref{prexfzxz}) is valid and the estimate in part (a) of Lemma~\ref{logexit1b} holds. 

Consider the Markov renewal process $(\xi^{ \varepsilon}_k, \tau^{\varepsilon}_k)$ on the state space 
$M = M_{r,\varkappa} = \Gamma_{r\varepsilon} \bigcup \Gamma_\varkappa$ with the starting point $x$. The process
is defined as follows. For $k = 0$,  $\xi^{\varepsilon}_0 = x$, 
$\tau^{\varepsilon}_0 = 0$. For $k \geq 1$, let
\[
\tilde{\tau}^{\varepsilon}_k = \inf\{t > \tau^{ \varepsilon}_{k-1}: X^{ \varepsilon}_t \in \Gamma_{2r \varepsilon} \}.
\] 
Then we define
\[
{\tau}^{\varepsilon}_k = \inf\{t > \tilde{\tau}^{\varepsilon}_{k}: X^{\varepsilon}_t \in M \},~~~\xi^{ \varepsilon}_k = 
X^{\varepsilon}_{{\tau}^{\varepsilon}_k}.
\] 
Define the random variable $K^{\varepsilon} $ as
\[
K^{\varepsilon} = \min\{k: \xi^{ \varepsilon}_k \in \Gamma_\varkappa \}.
\]
Thus 
\begin{equation} \label{oiop1}
\mathrm{E}_x \tau^{\varepsilon} (\Gamma_\varkappa) = \mathrm{E}_x  {\tau}^{\varepsilon}_{K^{\varepsilon} } = \sum_{k=0}^{\infty}
\mathrm{E}_x  \left((\tau^{\varepsilon}_{k+1} - \tau^{\varepsilon}_{k})\chi_{\{K^{\varepsilon} > k\}}\right) \leq 
\sup_{x \in \Gamma_{r \varepsilon}} \mathrm{E}_x \tau^{ \varepsilon}_1\sum_{k=0}^{\infty}
\mathrm{P}_x(K^{\varepsilon} > k) ,  
\end{equation}
where the inequality is obtained by conditioning on the $\sigma$-algebra of events determined by the time $\tau^{\varepsilon}_{k}$.
It follows from 
Lemma~\ref{epsde} (which bounds the expected transition time  from $\Gamma_{r\varepsilon}$ to $\Gamma_{2 r \varepsilon}$), part (a) of Lemma~\ref{logexit1b} (which bounds the expected time required to reach $\Gamma_{r\varepsilon} \bigcup \Gamma_\varkappa$ 
from $\Gamma_{2r \varepsilon}$), and the strong Markov property that
$\sup_{x \in \Gamma_{r \varepsilon}} \mathrm{E}_x \tau^{\varepsilon}_1$ is bounded from above by a constant. The sum in the right-hand side of (\ref{oiop1})
is bounded from above by $c({\varkappa}/{\varepsilon})^\gamma$ for some $c > 0$, as follows from the lower bound in (\ref{prexfzxz}).
We thus obtain (\ref{time1}). 

Next, we prove (\ref{bfb}). Consider sufficiently small $\varkappa$ so that
the upper bound in (\ref{prexfzxz}) holds with $\varkappa_2 = \varkappa$.
Pick an arbitrary $\delta > 0$. From the upper bound in (\ref{prexfzxz}), it follows that there exist $c >0$ and $\varepsilon_0$ such that,
for $0 < \varepsilon \leq \varepsilon_0$, 
\begin{equation} \label{emnn}
\sup_{x \in V_{0,r\varepsilon}} \mathrm{P}_x(K^{\varepsilon} \leq [c({\varkappa}/{\varepsilon})^\gamma]+1) \leq {\delta}/{2}. 
\end{equation}
Consider the generator of the process $X^\varepsilon_t$ in $V_{0, 2r \varepsilon}$. After the change of variables 
$\tilde{z} = z/\varepsilon$, its coefficients remain bounded, as follows from Lemma~\ref {geneNN}. This implies that 
we can pick $s' >0 $ such that  
\begin{equation} \label{juk}
\mathrm{P}_x(\tau^{\varepsilon}_1 \leq s') \leq 1/2~~{\rm for}~ x \in \Gamma_{r\varepsilon}.
\end{equation}
Consider $s < c s'\varkappa^\gamma/5$. Define $n(\varepsilon) = [c({\varkappa}/{\varepsilon})^\gamma]+1$. Let $E^{\varepsilon}$ be the event that $K^{\varepsilon} \geq n(\varepsilon)$
and $F^{\varepsilon}$  be the event that at least $1/5$ of 
the random variables $\tau^{ \varepsilon}_{1} - \tau^{ \varepsilon}_{0}$,...,
$\tau^{ \varepsilon}_{n(\varepsilon)} - \tau^{ \varepsilon}_{n(\varepsilon) -1}$ exceed $s'$. Observe that 
$\varepsilon^\gamma   {\tau^{\varepsilon} (\Gamma_\varkappa)}  \geq s$ on $E^{\varepsilon} \bigcap F^{\varepsilon}$. By (\ref{emnn}), the strong Markov property,
and (\ref{juk}),
\[
\mathrm{P}_x(E^{\varepsilon} \bigcap F^{\varepsilon} ) \geq \mathrm{P}_x(E^{\varepsilon}) - \mathrm{P}_x(E^{\varepsilon} \setminus F^{\varepsilon}) \geq 1 - \frac{\delta}{2} 
 - (\frac{1}{2})^{[\frac{4n(\varepsilon)}{5}]} C^{n(\varepsilon)}_{ [\frac{4n(\varepsilon)}{5}]}.
\]
Tne last term on the right-hand side can be made smaller than $\delta/2$ by selecting sufficiently small $\varepsilon$, which completes the
proof of the lemma. \qed
\\

From Lemma~\ref{pratfzzz} (which implies that the process starting near $\mathcal{M}$ reaches $V_{0,r\varepsilon}$ prior to reaching $\Gamma_\varkappa$ with probability close to one if $r$ is sufficiently large), part (b) of Lemma~\ref{mltime} (which estimates, from below, the time it takes to leave a small neighborhood of an attracting surface if we start in  $V_{0,r\varepsilon}$, and Lemma~\ref{proxnnn2}, we obtain
the following improved version of Lemma~\ref{proxnnn2} for attracting surfaces.
\begin{lemma} \label{imprv}
Let $\gamma > 0$. If $1 \ll t(\varepsilon) \ll \varepsilon^{-\gamma}$ and $f \in C_b(\mathbb{R}^d)$, then
\[
\lim_{ \varepsilon, \varkappa \downarrow 0} \sup_{x \in V_{0,\varkappa}} |\mathrm{E}_x f(X^\varepsilon_t) - \int_{\mathcal{M}} {f} d\pi_{\mathcal{M}}|  = 0.
\]
\end{lemma}

\section{Distribution of the process as it exits a neighborhood of an invariant surface} \label{dpen}

For small $\varkappa > 0$, define $\overline{\mu}^x_\varkappa$ to be the measure induced on $\Gamma_\varkappa$ by $X_{\tau(\Gamma_\varkappa)}$, where the process $X_t$ starts at $x \in V_{0, \varkappa}$. Observe that $\mathrm{P}_x(\tau(\Gamma_\varkappa) < \infty) < 1$  if $\gamma > 0$ and $x \in \Gamma_{\varkappa'}$ with $\varkappa' < \varkappa$, and therefore $\overline{\mu}^x_\varkappa$ is not a probability measure
in this case. For $x \notin \mathcal{M}$, define the probability measure ${\mu}^x_\varkappa$ as the normalized version of $\overline{\mu}^x_\varkappa$. The measure $\mu^{x,\varepsilon}_\varkappa$ is defined in the same way using the process $X^\varepsilon_t$ instead of $X_t$ (no normalization is needed in this case since the perturbed process reaches $\Gamma_\varkappa$ with probability one even in the case of an attracting surface). The following lemma holds for both attracting and repelling surfaces.
\begin{lemma} \label{exds} Suppose that ${\rm dim}(\mathcal{M}) < d-1$. 
For each $\varkappa > 0$ that is sufficiently small, 
there is a probability measure $\nu_\varkappa$ on $\Gamma_\varkappa$ such that for each $f \in C(\Gamma_\varkappa)$,
\begin{equation} \label{upr}
\lim_{\substack{{\rm dist} (x, \mathcal{M}) \rightarrow 0 \\ x \notin \mathcal{M}}} \int_{\Gamma_\varkappa} f d \mu^x_\varkappa = \int_{\Gamma_\varkappa} f d \nu_\varkappa,
\end{equation} 
\begin{equation} \label{upr2}
\lim_{\substack{{\rm dist} (x, \mathcal{M}) \rightarrow 0 \\ \varepsilon \downarrow 0}} \int_{\Gamma_\varkappa} f d \mu^{x,\varepsilon}_\varkappa = \int_{\Gamma_\varkappa} f d \nu_\varkappa.
\end{equation}
\end{lemma} 
\noindent
{\bf Remark.} It is important (and this is used in Section~\ref{mdis}) that the limit in (\ref{upr2}) is a joint limit and thus the measure
$\nu_\varkappa$ in the right-hand side does not depend on $\varepsilon$.
\\

The rest of this section is mostly devoted to the proof of Lemma~\ref{exds}. We will need the following result. 
\begin{lemma} \label{exle1} Let  $U' \subseteq U$ be open connected domains (in $\mathbb{R}^d$ or in a smooth $d$-dimensional surface embedded in a Euclidean space). Let  $\partial U = \partial U' = B$ be a smooth $(d-1)$-dimensional compact surface. Let $K \subset U'$   be a compact set. Let $H_t$ be a diffusion process on $U \bigcup B$ whose generator is uniformly elliptic on $U' \bigcup B$.
Let $\tau = \inf\{t: H_t \in B \}$ (this stopping time may be equal to infinity with positive probability). Let $\mu^x$ the the sub-probability measure on $B$ induced by $H_{\tau}$ starting at $x$, and let $p^x$ be its density with respect to the Lebesgue measure on $B$. Then there is a constant $c >0$ such  
that 
\begin{equation} \label{densiest}
p^x(\tilde{x}) \geq c,~~x \in K,~\tilde{x} \in B. 
\end{equation}
The constant $c$ 
can be chosen  in such a way that (\ref{densiest})
holds for all the operators that have the same ellipticity constant and bound on the $C$-norm of the coefficients in $U'$.
\end{lemma}
The estimate (\ref{densiest}) is likely known in the PDE literature but, since we could not find a proof, we sketched a version of this
estimate in \cite{FK21}. 
\\

Lemma~\ref{exle1} will be applied to the process $X_t$ but written in different coordinates. 
Namely, let us write the process $X_t$ in $(y,z)$ coordinates as $(Y_t, Z_t)$. We will also consider the process
$(\mathcal{Y}_t, \mathcal{Z}_t) = (Y_t, \ln(Z_t))$ with $\mathbf{x} = (y, \mathbf{z}) = (y, \ln(z))$ as an initial point. The process
can be considered in the domain 
\[
\mathcal{V}_n = \{(y,\mathbf{z}):  \frac{1}{\gamma} \ln(\varphi(y)) +  \mathbf{z}  \leq  n\} \subset \mathbf{S} \times \mathbb{R}.
\]
By Lemma~\ref{geneNN}, the generator of the process $(\mathcal{Y}_t, \mathcal{Z}_t) $ in  $(y,\mathbf{z})$ coordinates can be written as:
\[
\mathbf{L} u  = L_y u  + \alpha(y) \frac{\partial^2 u}{\partial \mathbf{z}^2} +   \beta(y) \frac{\partial u}{\partial \mathbf{z}}  + 
{\mathcal{D}}_y \frac{\partial u}{\partial \mathbf{z}} + e^{\mathbf{z}} (  {\mathcal{K}}_y u   +  \mathcal{N}_y 
\frac{\partial u}{\partial \mathbf{z}}  +  \sigma(y,\mathbf{z}) \frac{\partial^2 u}{\partial \mathbf{z}^2}) =: \mathbf{M} + e^{\mathbf{z}}
\mathbf{R}.
\]
Observe that the operator $\mathbf{M}$ is translation-invariant in $\mathbf{z}$, while the domains $\mathcal{V}_{n}$ with different 
$n$ can be obtained from each other via translations in the $\mathbf{z}$ coordinate.

Let us take $\mathbf{x} = (y, \mathbf{z})$ such that $\frac{1}{\gamma} \ln(\varphi(y)) +  \mathbf{z} = n-1$ as an initial point for
the process $(\mathcal{Y}_t, \mathcal{Z}_t)$, which takes values in the domain 
$\mathcal{V}_{n}$. 
In order to get a domain that does not depend on $n$, we can consider the process $(\mathcal{Y}_t, \mathcal{Z}_t - n)$ in the domain $\mathcal{V}_{0}$, which is now independent of $n$. The initial point $\mathbf{x}' = (y, \mathbf{z}') = (y, \mathbf{z}-n)$ now satisfies  
$\frac{1}{\gamma} \ln(\varphi(y)) +  \mathbf{z}' = -1$. The generator of the latter process is
$\mathbf{M} + e^{\mathbf{z}+n} \mathbf{R}$, which a small perturbation of the operator  $\mathbf{M}$ if $-n$ is sufficiently large,
and thus  Lemma~\ref{exle1} applies.

The result of Lemma~\ref{exle1} can be now interpreted in terms of the process $(\mathcal{Y}_t, \mathcal{Z}_t)$ and, consequently, in terms of the process $X_t = (Y_t, Z_t)$ on the domain $V_{0,e^n}$. Namely, Lemma~\ref{exle1} implies the following.
\begin{lemma} \label{coxxx} There exist  $c > 0$ and $N \in \mathbb{Z}$ such that, for each $n \leq N$,
the measures induced by the random variables $X_{\tau(\Gamma_{e^n})}$ with $X_0 = x \in \Gamma_{e^{n-1}}$, have a common (for all $x$) component 
$\mu_n$ on $\Gamma_{e^n}$ such that $\mu_n(\Gamma_{e^n}) \geq c$.
%
\end{lemma}

\noindent
{\it Proof of Lemma~\ref{exds}.} For a probability measure $\psi$
on $V_{0,\varkappa}$, define $\overline{\mu}^\psi_\varkappa$ to be the measure induced by $X_{\tau(\Gamma_\varkappa)}$, where $\psi$ serves as the initial distribution of the process $X_t$. (Thus $\overline{\mu}^{\delta_x}_\varkappa$, 
$x \in V_{0, \varkappa}$, 
coincides with $\overline{\mu}^x_\varkappa$ introduced above.) Similarly, for $\psi$ that is not concentrated on $\mathcal{M}$, 
let ${\mu}^\psi_\varkappa$ be
the normalized version of  $\overline{\mu}^\psi_\varkappa$. Also, we define the measure $\mu^{\psi, \varepsilon}_\varkappa$ in the same way,
using the process $X^\varepsilon_t$.

We start by proving (\ref{upr}) in the case when $\gamma < 0$. 
Take $N \in \mathbb{Z}$ such that 
$e^N \leq \varkappa$ and the result of Lemma~\ref{coxxx} holds. Consider $k \geq 1$, and suppose that $\psi$ is a probability distribution
on $\Gamma_{e^{N-k}}$. 
Since $\gamma < 0$, the measure $\overline{\mu}^\psi_{e^{N-k+1}}$ induced by the stopped process is a probability measure, i.e., 
${\mu}^\psi_{e^{N-k+1}} = \overline{\mu}^\psi_{e^{N-k+1}}$. Lemma~\ref{coxxx} provides a bound on the total variation distance between the induced measure corresponding to two initial measures:
\begin{equation} \label{proxab}
\rho({\mu}^{\psi_1}_{e^{N-k+1}}, {\mu}^{\psi_2}_{e^{N-k+1}}) \leq 1-c,~~~\psi_1, \psi_2 - {\rm measures}~{\rm on}~\Gamma_{e^{N-k}}.
\end{equation}
Using the strong Markov property, we obtain (if $k \geq 2$) that
\[
\rho({\mu}^{\psi_1}_{e^{N-k+2}}, {\mu}^{\psi_2}_{e^{N-k+2}}) \leq (1-c)^2,~~~\psi_1, \psi_2  - {\rm measures}~{\rm on}~\Gamma_{e^{N-k}},
\]
and, continuing by induction, we get
\[
\rho({\mu}^{\psi_1}_{\varkappa}, {\mu}^{\psi_2}_{\varkappa}) \leq
\rho({\mu}^{\psi_1}_{e^{N}}, {\mu}^{\psi_2}_{e^{N}}) \leq (1-c)^k,~~~\psi_1, \psi_2  - {\rm measures}~{\rm on}~\Gamma_{e^{N-k}}.
\]
Applying this to $\psi_1 = \mu^{x_1}_{e^{N-k}}$, $\psi_2 = \mu^{x_2}_{e^{N-k}}$, with $x_1, x_2 \in V_{0, e^{N-k}} \setminus \mathcal{M}$, we 
obtain
\[
\rho({\mu}^{x_1}_{\varkappa}, {\mu}^{x_2}_{\varkappa}) \leq
(1-c)^k,~~~x_1, x_2 \in V_{0, e^{N-k}} \setminus \mathcal{M}.
\]
Thus there is a limit, in total variation, for the measures $\mu^x_\varkappa$ as ${\rm dist} (x, \mathcal{M}) \rightarrow 0$, 
$x \notin \mathcal{M}$. The limit will be denoted by $\nu_\varkappa$. This justifies (\ref{upr}) in the case when $\gamma < 0$.

Next, we sketch the proof of (\ref{upr}) in the case when $\gamma > 0$. Here, the additional difficulty is that $\overline{\mu}^\psi_{e^{N-k+1}}$ is not
a probability measure for $\psi$ that is a measure  on $\Gamma_{e^{N-k}}$, and the required normalization may depend on $\psi$.
In order to circumvent this problem, observe that from Lemma~\ref{pratfzzz} it follows that 
$\overline{\mu}^{\psi_1}_{e^{N-k+1}}(\Gamma_{e^{N-k+1}})/\overline{\mu}^{\psi_2}_{e^{N-k+1}}(\Gamma_{e^{N-k+1}}) \rightarrow 1$ as
$N \rightarrow -\infty$ uniformly in $\psi_1, \psi_2$ $k \geq 1$, i.e., the normalization required is almost the same if $N$ is sufficiently close to $-\infty$.   From here,   using the definition of ${\mu}^{\psi_1}_{e^{N-k+1}}$, ${\mu}^{\psi_2}_{e^{N-k+1}}$  as the normalized measures,  it follows that (\ref{proxab}) still holds for some $c > 0$, provided that  $-N$ is sufficiently large. Using the Markov property and induction, it is then seen that, for each $\delta > 0$, there is $k > 0$ such that 
\[
\rho({\mu}^{x_1}_{e^{N}}, {\mu}^{x_2}_{e^{N}}) \leq \delta,~~~x_1, x_2 \in V_{0, e^{N-k}} \setminus \mathcal{M}
\]
if $-N$ is sufficiently large.
 From Lemma~\ref{pratf}, it follows that 
$\overline{\mu}^{\psi_1}_{\varkappa}(\Gamma_{\varkappa})/\overline{\mu}^{\psi_2}_{\varkappa}(\Gamma_{\varkappa}) \in [1/2, 2]$  if $\psi_1, \psi_2$ are measures on $\Gamma_{e^N}$ if $\varkappa$ is sufficiently small and $e^N \leq \varkappa$, which implies that the proximity
of $\psi_1$ and $\psi_2$ leads to the proximity of $\mu^{\psi_1}_\varkappa$ and $\mu^{\psi_2}_\varkappa$. 
Thus, for each $\delta > 0$, there is $k > 0$ such that 
\[
\rho({\mu}^{x_1}_{\varkappa}, {\mu}^{x_2}_{\varkappa}) \leq 2 \delta,~~~x_1, x_2 \in V_{0, e^{N-k}} \setminus \mathcal{M},
\]
if $-N$ is sufficiently large.
Thus, again, justifies the existence of a limit (denoted by $\nu_\varkappa$), in total variation, for the measures 
$\mu^x_\varkappa$ as ${\rm dist} (x, \mathcal{M}) \rightarrow 0$. This proves (\ref{upr}) in the case when $\gamma > 0$.

Next, let us sketch the proof of (\ref{upr2}). Given 
$\varkappa > 0$ for which (\ref{upr}) works and an arbitrary $\delta >0$, 
we can find $N  \in \mathbb{Z}$ with  $e^{N} < \varkappa$ such that
\[
|\int_{\Gamma_\varkappa} f d  {\mu}^x_\varkappa  - \int_{\Gamma_\varkappa} f d \nu_\varkappa| < \delta,~~x \in \Gamma_{N}.
\]
For $0 < \varkappa' < \varkappa$ and $x \in V_{\varkappa',\varkappa}$, define $\overline{\mu}^x_{\varkappa',\varkappa}$ to be the measure induced on  $\Gamma_\varkappa$ by the stopped process $ X_{\tau(\Gamma_\varkappa)}$ restricted to the event 
$\{\tau(\Gamma_\varkappa) < \tau(\Gamma_{\varkappa'}')\}$, where the process $X_t$ starts at $x$. Observe that  
$\overline{\mu}^x_{\varkappa',\varkappa}$ is not a probability measure, but we can normalize it to define the probability measure ${\mu}^x_{\varkappa',\varkappa}$. Also note that ${\mu}^x_{\varkappa',\varkappa}$ tends to ${\mu}^x_{\varkappa}$ as $\varkappa' \downarrow 0$, and therefore we can find a sufficiently small $\varkappa'$ such that
\[
|\int_{\Gamma_\varkappa} f d {\mu}^x_\varkappa  - \int_{\Gamma_\varkappa} f d {\mu}^x_{\varkappa',\varkappa}| < 
\delta,~~x \in \Gamma_{N},
\]
and therefore
\[
|\int_{\Gamma_\varkappa} f d {\mu}^x_{\varkappa',\varkappa} - \int_{\Gamma_\varkappa} f d \nu_\varkappa| < 2\delta,~~x \in \Gamma_{N}.
\]
From here, due to the proximity of the processes $X_t$ and $X^\varepsilon_t$ on finite time intervals, considering successive 
excursions of $X^\varepsilon_t$ between $\Gamma_N$ and $\Gamma_{\varkappa'}$ prior to reaching $\Gamma_\varkappa$, we obtain that 
\[
|\int_{\Gamma_\varkappa} f d {\mu}^{x,\varepsilon}_{\varkappa} - \int_{\Gamma_\varkappa} f d \nu_\varkappa| < 3\delta,~~x \in \Gamma_{N},
\]
provided that $\varepsilon$ is sufficiently small. 
This easily implies (\ref{upr2}).   \qed

\section{Metastable distributions for the perturbed process} \label{mdis}
In this section, we describe the distribution of the process $X^{x,\varepsilon}_t$ at different time scales.

\subsection{Asymptotics of transition probabilities between different attracting surfaces} \label{atpr}

First, we examine transition probabilities between small neighborhoods
of different attracting surfaces for the unperturbed process. Recall that $\overline{m}$ is the number of attracting surfaces, i.e., we have $\gamma_1 \geq ... \geq \gamma_{\overline{m}} > 0 > \gamma_{\overline{m}+1} \geq  ... \geq \gamma_m$. Here, we assume that $\overline{m} \geq 2$. 
The quantities $p^x_k = \mathrm{P}_x(E_k)$ defined in Section~\ref{ldup} clearly satisfy $\lim_{{\rm dist}(x, \mathcal{M}_i) \downarrow 0} p^x_j = 0$ if $0 \leq i,j \leq
\overline{m}$, $i \neq j$. Of interest, however, is the behavior of conditional probabilities $\mathrm{P}_x(E_j|E_i^c)$ as 
${\rm dist}(x, \mathcal{M}_i) \downarrow 0$, $x \notin  \mathcal{M}_i$, where $E_i^c$ is the complement of the event $E_i$. We have the following lemma.
\begin{lemma} \label{conlemj}  Suppose that $d - d_k > 1$ for each $k =1,...,m$. 
For $0 \leq i,j \leq
\overline{m}$, $i \neq j$, there are limits
\begin{equation} \label{conpr1}
q_{ij} = \lim_{\substack{{\rm dist} (x, \mathcal{M}_i) \rightarrow 0 \\ x \notin \mathcal{M}_i}}  \frac{p^x_j}{\sum_{k \neq i} p^x_k}.
\end{equation}
\end{lemma}
\proof Let the surfaces $\Gamma^k_\varkappa$ and the sets $V^k_{\varkappa_1, \varkappa_2}$ be defined as in Section~\ref{tpni} but with the surface $\mathcal{M}_k$ instead of a generic surface $\mathcal{M}$. Let $\varkappa$ be sufficiently small so that 
(\ref{upr}) holds for $\mathcal{M}_i$.  For $x \in \Gamma^i_\varkappa$, define
\[
f(x) = \mathrm{P}_x(E_i) = \mathrm{P}_x(\lim_{t \rightarrow \infty}  
{\rm dist}(X_t, \mathcal{M}_i) = 0),~~~g(x) = \mathrm{P}_x(E_j) = \mathrm{P}_x(\lim_{t \rightarrow \infty}  
{\rm dist}(X_t, \mathcal{M}_j) = 0).
\]
Due to non-degeneracy of the process $X_t$ in $D$ and Lemma~\ref{cohhi}, the functions 
$f$ and $g$ take values in $(0,1)$ and are continuous on $\Gamma^i_\varkappa$. By the strong Markov property,
\[
\frac{p^x_j}{\sum_{k \neq i} p^x_k} = \frac{ \mathrm{P}_x(E_j)}{\mathrm{P}(E_i^c)} = 
\frac{\int_{\Gamma^i_\varkappa} g  d\overline{\mu}^x_{i, \varkappa}}{\int_{\Gamma^i_\varkappa} (1-f)  
d \overline{\mu}^x_{i, \varkappa}} = \frac{\int_{\Gamma^i_\varkappa} g  d{\mu}^x_{i, \varkappa}}{\int_{\Gamma^i_\varkappa} (1-f)  
d {\mu}^x_{i, \varkappa}},
\]
where $\overline{\mu}^x_{i, \varkappa}$ and ${\mu}^x_{i, \varkappa}$  are defined as in Section~\ref{dpen} but for the specific surface $\mathcal{M}_i$. By Lemma~\ref{exds}, 
\[
\lim_{\substack{{\rm dist} (x, \mathcal{M}_i) \rightarrow 0 \\ x \notin \mathcal{M}_i}}  \frac{p^x_j}{\sum_{k \neq i} p^x_k} =
\lim_{\substack{{\rm dist} (x, \mathcal{M}_i) \rightarrow 0 \\ x \notin \mathcal{M}_i}}   \frac{\int_{\Gamma^i_\varkappa} g  
d{\mu}^x_{i, \varkappa}}{\int_{\Gamma^i_\varkappa} (1-f)  
d {\mu}^x_{i, \varkappa}} = \frac{\int_{\Gamma^i_\varkappa} g  
d{\nu}_{i, \varkappa}}{\int_{\Gamma^i_\varkappa} (1-f)  
d {\nu}_{i, \varkappa}} =: q_{ij}.
\]
Observe that the last ratio does not depend on $\varkappa$ since the left-hand side does not depend on $\varkappa$. This completes the proof of the lemma.
  \qed
\\
\\
{\bf Remark.} The quantity  ${p^x_j}/{\sum_{k \neq i} p^x_k}$ in the right-hand side of (\ref{conpr1}) is the probability that the process
$X_t$ starting at $x$ and conditioned on not going to $\mathcal{M}_i$ in the limit goes to $\mathcal{M}_j$ in the limit. The  process, conditioned on not going to $\mathcal{M}_i$, is also a diffusion process. Its generator can be obtained from $L$ 
via the Doob transform. Namely, it is possible to show that there is a unique bounded function $h^i \in C^2(D)$ that satisfies
\[
L h^i(x) = 0,~~x \in D;~~\lim_{{\rm dist}(x, \mathcal{M}_i) \downarrow 0} h^i(x) = 1;~~
\lim_{{\rm dist}(x, \mathcal{M}_j) \downarrow 0} h^i(x) = 0,~~1 \leq j \leq \overline{m},~j \neq i.
\]
The generator of the conditioned process is the operator $L^i$, where $L^i u = L(h^i u)/h^i$. The quantities $q_{ij}$ can be
obtained by solving, in the class of bounded $C^2(D)$ functions, the equation
\[
L^i u_{ij}(x) = 0,~~x \in D;~~ \lim_{{\rm dist}(x, \mathcal{M}_j) \downarrow 0} u_{ij} = 1;~~ 
\lim_{{\rm dist}(x, \mathcal{M}_k) \downarrow 0} u_{ij}(x) = 0,~~1 \leq k \leq \overline{m},~k \neq i,j,
\]
and then taking the limit
\[
q_{ij} = \lim_{{\rm dist}(x, \mathcal{M}_i) \downarrow 0} u_{ij}(x).
\]
\\

Next, let us explore transitions between small neighborhoods of the attracting surfaces for the perturbed process.
Lemma~\ref{pratfzzz} can be applied to each of the surfaces $\mathcal{M}_k$, $1 \leq k \leq \overline{m}$. Let $r$ be sufficiently large and $\varkappa > 0$ be sufficiently small so that, for each $1 \leq k \leq \overline{m}$, the second
inequality in (\ref{prexfzxz}) holds with $\eta = 1$ and $\varkappa_2 = \varkappa$. In particular, taking $\varkappa_1 = r \varepsilon$,
we get 
\begin{equation} \label{prexfzxz22}
\mathrm{P}_x 
( X^{\varepsilon}_{ \tau^{\varepsilon} (\Gamma^k_{r \varepsilon} \bigcup \Gamma^k_{{\varkappa}})} 
\in \Gamma^k_{{\varkappa}}) \leq 
\frac{2 \zeta^\gamma - (r \varepsilon)^\gamma}{{\varkappa}^\gamma - {(r \varepsilon)}^\gamma}
\end{equation}  
for all sufficiently small $\varepsilon$, 
provided that $x \in \Gamma^k_\zeta$ with $ r \varepsilon \leq \zeta \leq {\varkappa}$.

Define $\overline{\Gamma}_{r \varepsilon} = \bigcup_{k=1}^{\overline{m}} \Gamma^k_{r \varepsilon}$ and 
$\overline{\Gamma}_{r \varepsilon}^i =  \overline{\Gamma}_{r \varepsilon} \setminus  {\Gamma}_{r \varepsilon}^i$.

\begin{lemma} \label{trapl} Suppose that $d - d_k > 1$ for each $k =1,...,m$.  For the constants $q_{i j} > 0$, $1 \leq i, j \leq \overline{m}$, $i \neq j$, defined in Lemma~\ref{conlemj}, we have
\[
\lim_{\varepsilon \downarrow 0} \mathrm{P}_x(X^{ \varepsilon}_{\tau^\varepsilon(\overline{\Gamma}_{r \varepsilon}^i)} \in 
\Gamma^j_{r \varepsilon}) = q_{i j},~~x \in \Gamma^i_{r \varepsilon}.
\]
The convergence is uniform in $x \in \Gamma^i_{r \varepsilon}$. 
\end{lemma} 
\proof Let $\varkappa$ be sufficiently small so that 
(\ref{upr}), (\ref{upr2}) hold for $\mathcal{M}_i$. Let $A_n^\varepsilon$, $n \geq 1$, be the event that the process $X^\varepsilon_t$ makes
at least $n$ transitions from $ \Gamma^i_{r \varepsilon}$ to $\Gamma^i_{\varkappa}$ prior to reaching 
$\overline{\Gamma}_\varepsilon \setminus \Gamma^i_{r \varepsilon}$.  For $x \in \Gamma^i_\varkappa$, define
\[
f^\varepsilon(x) = \mathrm{P}_x(X^\varepsilon_{\tau^\varepsilon(\overline{\Gamma}_{r\varepsilon})} \in \Gamma^i_{r \varepsilon}),~~~g^\varepsilon(x) = 
\mathrm{P}_x(X^\varepsilon_{\tau^\varepsilon(\overline{\Gamma}_{r\varepsilon})} \in \Gamma^j_{r \varepsilon}).
\]
Let ${\mu}^{x,\varepsilon, n}_{i, \varkappa}$,   be the measure on $\Gamma^i_{\varkappa}$ induced 
by the process $X^\varepsilon_t$ (considered on the event $A^\varepsilon_n$) stopped after the $n$-th transition from $\Gamma^i_{r \varepsilon}$ to $\Gamma^i_\varkappa$, 
with the starting point $x \in \Gamma^i_{r \varepsilon}$. Thus ${\mu}^{x,\varepsilon, 1}_{i, \varkappa}$
is  the same as the measure ${\mu}^{x,\varepsilon}_{\varkappa}$ defined as in Section~\ref{dpen}, but now for the specific surface $\mathcal{M}_i$. For $n \geq 2$, ${\mu}^{x,\varepsilon, n}_{i, \varkappa}$ are sub-probability measures since $\mathrm{P}_x (A^\varepsilon_n) <1$ for $n \geq 2$. By the strong Markov property,
\begin{equation} \label{eqoi}
 \mathrm{P}_x(X^{ \varepsilon}_{\tau^\varepsilon(\overline{\Gamma}_{r \varepsilon}^i)} \in 
\Gamma^j_{r \varepsilon}) = \sum_{n =1}^\infty \int_{\Gamma^i_\varkappa} g^\varepsilon d {\mu}^{x,\varepsilon, n}_{i, \varkappa}.
\end{equation}
From non-degeneracy of $X_t$ in D, the proximity of $X_t$ and $X^\varepsilon_t$ on finite time intervals, Lemmas~\ref{pratf}, 
\ref{pratfzzz}, and  \ref{cohhi}, it follows that  
$\lim f^\varepsilon (x) = f(x)$, $\lim g^\varepsilon (x) = g(x)$ uniformly in $x \in \Gamma^i_\varkappa$, where $f$ and $g$ were defined in the proof of Lemma~\ref{conlemj}. By Lemma~\ref{exds},
the measures ${\mu}^{x,\varepsilon, 1}_{i, \varkappa}$ converge, as $\varepsilon \downarrow 0$, $x \in \Gamma^i_{r \varepsilon}$, to the measure $\nu_{i, \varkappa}$. Therefore,
\[
\lim_{ \varepsilon \downarrow 0} \mathrm{P}_x (A^\varepsilon_2) =
\lim_{ \varepsilon \downarrow 0}  \int_{\Gamma^i_\varkappa} f^\varepsilon d {\mu}^{x,\varepsilon, 1}_{i, \varkappa} = \int_{\Gamma^i_\varkappa} f d \nu_{i, \varkappa} =:r < 1,
\]
where the limit is uniform in $x \in \Gamma^i_{r \varepsilon}$. 
Therefore,
from the strong Markov property and Lemma~\ref{exds}, it follows that the measures ${\mu}^{x,\varepsilon, 2}_{i, \varkappa}$, $x \in 
\Gamma^i_{r \varepsilon}$, converge, as  $\varepsilon \downarrow 0$, to the sub-probability measure $r \nu_{i, \varkappa}$. Continuing by induction, we see
that the measures ${\mu}^{x,\varepsilon, n}_{i, \varkappa}$, $x \in \Gamma^i_{r \varepsilon}$, converge, as  $\varepsilon \downarrow 0$, 
to  $r^{n-1} \nu_{i, \varkappa}$. Therefore, since the terms in the right-hand side of (\ref{eqoi}) decay exponentially uniformly in $x \in \Gamma^i_{r \varepsilon}$, we obtain
\[
\lim_{ \varepsilon \downarrow 0} 
\mathrm{P}_x(X^{ \varepsilon}_{\tau^\varepsilon(\overline{\Gamma}_{r \varepsilon}^i)} \in 
\Gamma^j_{r \varepsilon}) = \sum_{n =1}^\infty r^{n-1} \int_{\Gamma^i_\varkappa} g d {\nu}_{i, \varkappa} =
\frac{\int_{\Gamma^i_\varkappa} g  
d{\nu}_{i, \varkappa}}{\int_{\Gamma^i_\varkappa} (1-f)  
d {\nu}_{i, \varkappa}}.
\]
This coincides with the expression for $q_{ij}$ obtained in the proof of Lemma~\ref{conlemj}.
\qed

\subsection{The case when at least one of the invariant surfaces is attracting}

Here, we assume that $\gamma_1 > 0$, i.e., $\mathcal{M}_1$ is attracting. Recall that $\overline{m}$ is such that
$\gamma_1 \geq ... \geq \gamma_{\overline{m}} > 0 > \gamma_{\overline{m}+1} \geq ... \geq \gamma_m$. Observe that the functions $p^x_k = \mathrm{P}_x(E_k)$, $1 \leq k \leq \overline{m}$, defined in Section~\ref{ldup} are discontinuous on $\mathcal{M}_{\overline{m}+1} \bigcup... \bigcup \mathcal{M}_m$  since
the repelling surfaces are invariant for the unperturbed process. However, from Lemma~\ref{exds} it follows that these functions can be
extended from $D \bigcup \mathcal{M}_1 \bigcup... \bigcup \mathcal{M}_{\overline{m}}$ to $\mathcal{M}_{\overline{m}+1} \bigcup... \bigcup \mathcal{M}_m$ by continuity. We will use this continuous version of $p^x_k$, $1 \leq k \leq \overline{m}$ below.
\begin{lemma} \label{mmle} Suppose that $d - d_k > 1$ for each $k =1,...,m$.

(a) If
$t(\varepsilon) \gg 1$, then, for  each compact 
$K \subset D \bigcup \mathcal{M}_1 \bigcup... \bigcup \mathcal{M}_{\overline{m}}$ and each $\varkappa > 0$,
\begin{equation} \label{mm2}
\lim_{\varepsilon \downarrow 0} \sup_{x \in K} 
\mathrm{P}_x(\min(\tau^\varepsilon(V^1_{0,\varkappa}),...,  \tau^\varepsilon(V^{\overline{m}}_{0,\varkappa})) > t(\varepsilon)) = 0.
\end{equation}

(b)
If  $t(\varepsilon) \gg |\ln \varepsilon|$, then, for each compact $K \subset \mathbb{R}^d$,
\begin{equation} \label{mm3}
\lim_{\varepsilon \downarrow 0} \sup_{x \in K} 
\mathrm{P}_x(\min(\tau^\varepsilon(V^1_{0,r \varepsilon}),...,  \tau^\varepsilon(V^{\overline{m}}_{0,r \varepsilon})) > t(\varepsilon)) = 0.
\end{equation}

(c)
Let $K \subset \mathbb{R}^d$ be compact. Then
\begin{equation} \label{mm1}
\lim_{\varepsilon \downarrow 0}  
\mathrm{P}_x (
\tau^\varepsilon(V^k_{0,r \varepsilon}) = \min(\tau^\varepsilon(V^1_{0,r \varepsilon}),...,  \tau^\varepsilon(V^{\overline{m}}_{0,r \varepsilon}))) = p^x_k~~~uniformly~~in~~x \in K.
\end{equation}
\end{lemma}
\proof
First, consider the unperturbed process $X_t$.
Recall from the proof of Theorem~\ref{unpnn} that we defined $p^x_{k,\varkappa} = \mathrm{P}_x(\tau(V^k_{0,\varkappa}) = 
\min(\tau(V^1_{0,\varkappa}),...,  \tau(V^{\overline{m}}_{0,\varkappa})))$. Observe that the process $X_t$ is non-degenerate in $D$, does not escape to infinity (assumption (d)),  and the time spent during each visit to a sufficiently small  
neighborhood of $\mathcal{M}_{\overline{m}+1}\bigcup ... \bigcup \mathcal{M}_m$ has bounded expectation (Lemma~\ref{ttlrs}).
Therefore, for each $\varkappa, \delta > 0$, there is $t > 0$ such that
\[
\mathrm{P}_x(\min(\tau(V^1_{0,\varkappa}),...,  \tau(V^{\overline{m}}_{0,\varkappa})) \leq t) \geq 1 -\delta,~~~x \in K.
\]
From the proximity of $X_t$ and $X^\varepsilon_t$ on finite time intervals, we conclude that, for each $\varkappa, \delta > 0$, there is $t > 0$ such that
\[
\mathrm{P}_x(\min(\tau^\varepsilon(V^1_{0,\varkappa}),...,  \tau^\varepsilon(V^{\overline{m}}_{0,\varkappa})) \leq t) 
\geq 1 -\delta,~~~x \in K,
\]
provided that $\varepsilon$ is sufficiently small. This justifies (\ref{mm2}). 
Next, (\ref{mm3}) follows from (\ref{mm2}), (\ref{prexfzxz22}), and part (a) of Lemma~\ref{logexit1b}. 

Due to continuity of $p^x_k$, it is sufficient to prove (\ref{mm1}) for
$K \subset D \bigcup \mathcal{M}_1 \bigcup... \bigcup \mathcal{M}_{\overline{m}}$.
 Recall from the proof of Theorem~\ref{unpnn} that $p^x_{k,\varkappa} \rightarrow p^x_k$
as $\varkappa \downarrow 0$ uniformly on $K$. From the proximity of $X_t$ and $X^\varepsilon_t$ on finite time intervals, it follows that
\[
\lim_{\varepsilon \downarrow 0}  
\mathrm{P}_x (
\tau^\varepsilon(V^k_{0,\varkappa}) = \min(\tau^\varepsilon(V^1_{0,\varkappa}),...,  \tau^\varepsilon(V^{\overline{m}}_{0,\varkappa}))) = p^x_{k,\varkappa}~~~{\rm uniformly}~~{\rm in}~~x \in K.
\]
The result in (\ref{mm1}) now follows once we observe  that, by (\ref{prexfzxz22}), if $X^\varepsilon_t$ starts in $V^k_{0,\varkappa}$ with small $\varkappa$, then
it reaches $V^k_{0, r \varepsilon}$ with probability close to one prior to visiting small neighborhoods of other attracting surfaces.
\qed
\\

Next, we  describe the distribution of the perturbed process  at time scales $t(\varepsilon)$ that satisfy $1 \ll t(\varepsilon) \ll 
\varepsilon^{-\gamma_{\overline{m}}}$.  Recall that $\pi_{\mathcal{M}_k}$ are the invariant measures for the process $X^x_t$ considered as a process on $\mathcal{M}_k$.  

\begin{theorem} \label{iuo} Suppose that $d - d_k > 1$ for each $k =1,...,m$. If $1 \ll t(\varepsilon) \ll \varepsilon^{-\gamma_{\overline{m}}}$ and $x \in D \bigcup \mathcal{M}_1 \bigcup... \bigcup \mathcal{M}_{\overline{m}}$, then the distribution of 
$X^{x,\varepsilon}_{t(\varepsilon)}$ converges to the measure $\sum_{k=1}^{\overline{m}} p^x_k \pi_{\mathcal{M}_k}$, where the coefficients
$p^x_k$ are defined in Section~\ref{ldup}. For $x \in \mathcal{M}_{\overline{m}+1} \bigcup ... \bigcup \mathcal{M}_m$, the same conclusion 
holds if $|\ln(\varepsilon)| \ll t(\varepsilon) \ll \varepsilon^{-\gamma_{\overline{m}}}$.
\end{theorem}
\proof From Lemmas~\ref{mmle}, \ref{pratfzzz}, and \ref{mltime} (part (b)), it follows that, for sufficiently small $\varkappa > 0$, 
\[
\lim_{\varepsilon \downarrow 0} \mathrm{P}_x( X^\varepsilon_{t(\varepsilon)/2} \in V^k_{0,\varkappa}) = p^x_k.
\] 
 Therefore, by Lemma~\ref{imprv}, for each $f \in C_b(\mathbb{R}^d)$,
\[
\lim_{ \varepsilon \downarrow 0} |\mathrm{E}_x f(X^\varepsilon_{t(\varepsilon)}) -\sum_{k=1}^{\overline{m}} p^x_k  \int_{\mathcal{M}} {f} d\pi_{\mathcal{M}_k}|  = 0.
\]
This justifies the claim of the theorem.  \be
\qed
\\

Next, let us explore the behavior of $X^{x,\varepsilon}_t$ at longer time scales. Assume that $\varepsilon^{-\gamma_{l+1}} \ll t(\varepsilon)
\ll \varepsilon^{-\gamma_{l}}$, where $l+1 \leq \overline{m}$ and $l \geq 1$ (this is possible if $\overline{m} \geq 2$).  
Consider the  discrete-time 
Markov chain $Z_n^x$ on $\{1,...,\overline{m}\}$ with transition probabilities $q_{ij}$ defined in Lemma~\ref{conlemj} (it follows from
Lemma~\ref{conlemj} that these do form a stochastic matrix). We take the vector $(p^x_1,...,p^x_{\overline{m}})$ as the initial distribution
of the Markov chain. Let $\tau_l$ be the hitting time of the set $\{1,...l\}$, and define
\[
p^{x,l}_k = \mathrm{P}(Z^x_{\tau_l} = k),~~~k =1,...,l.
\]
\begin{theorem} \label{btto} Suppose that $d - d_k > 1$ for each $k =1,...,m$. If $\varepsilon^{-\gamma_{l+1}} \ll t(\varepsilon)
\ll \varepsilon^{-\gamma_{l}}$ and $x \in \mathbb{R}^d$, then the distribution of 
$X^{x,\varepsilon}_{t(\varepsilon)}$ converges to the measure $\sum_{k=1}^{l} p^{x,l}_k \pi_{\mathcal{M}_k}$. 
\end{theorem}
\proof From Lemma~\ref{mmle} (part (c)), it follows that
\[
\mathrm{P}_x(X^\varepsilon_{\tau^\varepsilon(V^1_{0, r\varepsilon} \bigcup ... \bigcup V^{\overline{m}}_{0, r \varepsilon})} \in V^k_{0, r\varepsilon}) 
\rightarrow p^{x}_k,~~~{\rm as}~~\varepsilon \downarrow 0,~~1 \leq k \leq \overline{m}. 
\]
From Lemma~\ref{trapl}, it then follows that 
\[
\mathrm{P}_x(X^\varepsilon_{\tau^\varepsilon(V^1_{0, r\varepsilon} \bigcup ... \bigcup V^l_{0, r \varepsilon})} \in V^k_{0, r\varepsilon}) 
\rightarrow p^{x,l}_k,~~~{\rm as}~~\varepsilon \downarrow 0,~~1 \leq k \leq l. 
\]
The set $V^1_{0, r\varepsilon} \bigcup ... \bigcup V^l_{0, r \varepsilon}$ is reached after the first transition from $x$ to 
$V^1_{0, r\varepsilon} \bigcup ... \bigcup V^{\overline{m}}_{0, r \varepsilon}$, which takes time $o(\varepsilon^\delta)$ as long as
 $\delta < 0$, by Lemma~\ref{mmle} (part (b)), and a finite number of transitions between the neighborhoods of different attracting surfaces 
 , which take time $O(\varepsilon^{-\gamma_{l+1}})$ by Lemma~\ref{mltime}. Therefore,
\[
\mathrm{P}_x({\tau^\varepsilon(V^1_{0, r\varepsilon} \bigcup ... \bigcup V^l_{0, r \varepsilon})} \leq \frac{t(\varepsilon)}{2} ) 
\rightarrow 1~~~{\rm as}~~\varepsilon \downarrow 0. 
\]
From Lemma~\ref{mltime} (part (b)), it then follows that, for sufficiently small $\varkappa > 0$, 
\[
\lim_{\varepsilon \downarrow 0} \mathrm{P}_x( X^\varepsilon_{t(\varepsilon)/2} \in V^k_{0,\varkappa}) = p^x_k,~~~1 \leq k \leq l.
\] 
The result now follows from Lemma~\ref{imprv}.
\qed
\\

Next, we explore the longest time scales. Here, for simplicity, we assume that $\gamma_1 > \gamma_2$. 
\begin{theorem} \label{ltsc} Suppose that $d - d_k > 1$ for each $k =1,...,m$. If $\gamma_1 > \gamma_2$, 
$ t(\varepsilon) \gg \varepsilon^{-\gamma_{1}}$,
 and $x \in \mathbb{R}^d$, then the distribution of 
$X^{x,\varepsilon}_{t(\varepsilon)}$ converges to the measure $\pi_{\mathcal{M}_1}$. 
\end{theorem}
\proof Take an arbitrary function $t_1(\varepsilon)$ that satisfies $\varepsilon^{-\gamma'} \ll t_1(\varepsilon) \ll \varepsilon^{-\gamma_1}$
with some $\gamma'$ such that $\max(0,\gamma_2) < \gamma' < \gamma_1$. Let $\mu^{x,\varepsilon}$ be the distribution of $X^{x,\varepsilon}_{t(\varepsilon) - t_1(\varepsilon)}$. Then the distribution of  $X^{x,\varepsilon}_{t(\varepsilon)}$ agrees with that of $X^{\mu^{x,\varepsilon},\varepsilon}_{t_1(\varepsilon)}$ (the process whose initial distribution is $\mu^{x,\varepsilon}$ rather than concentrated in one point). 
The same proof as for Theorem~\ref{iuo} (if $\gamma_2 < 0$) or Theorem~\ref{btto} (if $\gamma_2 > 0$) applies; the only difference is that now we need to bound the time it takes the process starting with the distribution $\mu^{x,\varepsilon}$ to reach a small neighborhood of 
one of the attracting surfaces, i.e., it is sufficient to establish that 
\[
\mathrm{P}_{\mu^{x,\varepsilon}}({\tau^\varepsilon(V^1_{0, r\varepsilon} \bigcup ... \bigcup V^{\overline{m}}_{0, r \varepsilon})} \leq \frac{t_1(\varepsilon)}{2} ) 
\rightarrow 1~~~{\rm as}~~\varepsilon \downarrow 0. 
\]
The latter easily follows from Lemma~\ref{mmle}, part (b). 
\qed
\\
\\
{\bf Remark.} At the `transitional' time scales ($t(\varepsilon) \sim \varepsilon^{-\gamma_l}$ with some $l \leq \overline{m}$), the
limiting distribution will also be a linear combination of the measures $\pi_{\mathcal{M}_k}$, $k \leq l$. Specifying the coefficients in this linear 
combination is not difficult if one takes into account that the time 
it takes the process $X^{x,\varepsilon}_t$ to exit  a neighborhood of an attracting surface is nearly exponentially distributed (after proper normalization) when
$\varepsilon$ is small. The statement on the exponential distribution, however, requires a separate argument that is not provided here. 
A similar result for transitional time scales was proved in the case of randomly perturbed dynamical systems in \cite{LLS1}, \cite{LLS2}. 

\subsection{The case when all the  components of the boundary are repelling}
Here, we assume that $\gamma_1 < 0$, i.e., all the components of the boundary are repelling.  Recall from Section~\ref{ldup} that $\mu$ is the
invariant measure for the unperturbed process $X^x_t$ on $D$.  
\begin{theorem} \label{hujuh4}  Suppose that $d - d_k > 1$ for each $k =1,...,m$. If $t(\varepsilon) \rightarrow \infty$ as $\varepsilon \downarrow 0$, then
the distribution of $X^{\varepsilon}_{t(\varepsilon)}$ converges to $\mu$ uniformly with respect to initial point on a given compact set. Namely, for each compact $K \subset D$ and each $f \in C_b(D)$, 
\begin{equation} \label{prox2}
\lim_{\varepsilon \downarrow 0} \sup_{x \in K} |\mathrm{E}_x f(X^{\varepsilon}_{t(\varepsilon)}) - \int_D f d\mu| =0.
\end{equation}
Formula (\ref{prox2}) holds with $K \subset \mathbb{R}^d$ and $f \in C_b(\mathbb{R}^d)$ if  $t(\varepsilon) \gg |\ln(\varepsilon)|$. 
\end{theorem}  
\proof Consider, first, the case where   $K \subset D$.
Take two $(d-1)$-dimensional spheres $F$ and $G$ such that the closed ball bounded by $F$ is contained in $D$ and $G$ 
is contained in 
the interior of this ball. Define, inductively,
\[
\sigma^\varepsilon_0 = \inf \{t \geq 0: X^\varepsilon_t \in F \},~~~
\hat{\sigma}^\varepsilon_0 = \inf \{t \geq \sigma^\varepsilon_0: X^\varepsilon_t \in G \},
\]
\[
\sigma^\varepsilon_n = \inf \{t \geq \hat{\sigma}^\varepsilon_{n-1}: X^\varepsilon_t \in F \},~~~
 \hat{\sigma}^\varepsilon_n = \inf \{t \geq \sigma^\varepsilon_n: X^\varepsilon_t \in G \}.
\]
Define $\mathbf{X}^\varepsilon_n =
X^\varepsilon_{\sigma^\varepsilon_n}$, $n \geq 0$.
Thus, assuming that the initial point belongs to $F$, the process $(\mathbf{X}^\varepsilon_n, \sigma^\varepsilon_n)$, $n \geq 0$,
is a Markov renewal process on $F$. (We use the arrival times, rather than the inter-arrival times, in the notation for the renewal process.) We have the following bounds on the inter-arrival times: $\mathrm{E}_x 
(\sigma_1^\varepsilon)^2 \leq c_1$, $\mathrm{P}_x(\sigma_1^\varepsilon \leq c_2) \leq 1/2$ for some positive constants $c_1$ and $c_2$, as follows from Lemma~\ref{semom} and the non-degeneracy of the process $X_t$ in $D$. From these bounds it follows (see, e.g., \cite{IK})  that the following property holds: for each $\delta > 0$, there is $L > 0$ such that, for each  $a \geq 0$, 
\begin{equation} \label{freq}
\mathrm{P}_x (\{\sigma^\varepsilon_n \in [a, a+L]~~{\rm for}~{\rm some}~n \geq 0\}) \geq 1 - \delta,~~x \in F.  
\end{equation}
The process $X^\varepsilon_t$ starting at $x \in K$ arrives to $F$ before time $L$ with probability close to one if $L$ is large enough, and therefore (\ref{freq}) holds for $x \in K$. 

Fix $\delta > 0$. By Theorem~\ref{repth}, there is $L_1$ be such that
\[
\sup_{x \in K} |\mathrm{E}_x f(X_{t}) - \int_D f d\mu| \leq \delta
\]
for $t \geq L_1$. Then, from the proximity of $X^\varepsilon_t$ and $X_t$ on finite time intervals, it follows that 
\begin{equation} \label{2de}
\sup_{x \in K} \sup_{t \in [L_1, L_1+L ]} |\mathrm{E}_x f(X^\varepsilon_{t}) - \int_D f d\mu| \leq 2\delta
\end{equation}
for sufficiently small $\varepsilon$. Choose $a = t(\varepsilon)  - L_1  -L $ in (\ref{freq}). Then, from the strong Markov property, 
(\ref{freq}) (with $F$ replaced by $K$), and (\ref{2de}), we obtain 
\[
\sup_{x \in K}   |\mathrm{E}_x f(X^\varepsilon_{t(\varepsilon)}) - \int_D f d\mu| \leq 3\delta.
\]
Since $\delta > 0$ was arbitrary, this implies the first statement of the lemma. The second statement (with  $K \subset \mathbb{R}^d$ and  $t(\varepsilon) \gg |\ln(\varepsilon)|$) easily follow the first statement and the fact that it takes at most logarithmic time for the process to leave a small neighborhood of a repelling surface (Lemma~\ref{ttlrs}).  \qed 
\\

Observe that Theorem~\ref{mnt2}  follows from Theorems~\ref{iuo}, \ref{btto}, \ref{ltsc}, and \ref{hujuh4} via the representation of the solution to Cauchy problem (\ref{direq2}) in the form  $u^\varepsilon(t(\varepsilon), x) = 
\mathrm{E}_x g(X^\varepsilon_{t(\varepsilon)})$.  
\\
\\
{\bf Acknowledgements.} Both authors were supported by the NSF grant DMS-2307377. Leonid Koralov was supported by the Simons Foundation Grant MP-TSM-00002743. 
\\

\end{document}